\documentclass{amsart}[11pt]
\usepackage{geometry,amsmath,color}
\usepackage{khodesh}

\renewcommand{\L}{\mathcal{L}}
\usepackage{mathpazo}
\begin{document}

\title{A spectral sequence for Lagrangian Floer homology}
\author{Reza Rezazadegan}

\begin{abstract}

 We prove the existence of a spectral sequence for Lagrangian Floer homology which  converges to the  Floer homology of the image of a Lagrangian submanifold under multiple fibred Dehn twists. 
 The $E_1$ term of the sequence is given by the hypercube of ``resolutions'' of the Dehn twists involved. The proof relies on the exact triangle for fibered Dehn twists due to Wehrheim and Woodward.

As applications we obtain a spectral sequences from Khovanov homology to symplectic Khovanov homology.
 Also when a 3-manifold $M$ is given by gluing two handlebodies by a surface diffeomorphism $\phi$,
 we obtain a spectral sequence converging to the Heegaard-Floer homology of $M$, whose $E_1$  term is  a  hypercube obtained from different ways of resolving the Dehn twists  in $\phi$.  This latter sequence generalizes the spectral sequence of branched double covers to general closed 3-manifolds (i.e. those which are not branched double covers of links) however its $E_2$ term  is not a 3-manifold invariant. This gives upper bounds on the rank of HF-hat of closed 3-manifolds.  

 \end{abstract}

\maketitle

\section{Introduction}

It is a well-known fact of homological algebra that whenever we have an iterated mapping cone of chain maps between complexes $C_{i}$ then there is a spectral sequence whose $E_{1}$ term is given by the direct sum of the cohomologies of the $C_{i}$ and converges to the cohomology of the iterated mapping cone.  In this paper we apply this principle to mapping cones arising from fibered Dehn twists along coisotropic submanifolds of symplectic manifolds. Such Dehn twists arise naturally in a variety of contexts and give rise to actions of the braid and mapping class groups by symplectomorphisms and therefore to link and 3-manifold invariants. 
 
\begin{figure}
\includegraphics{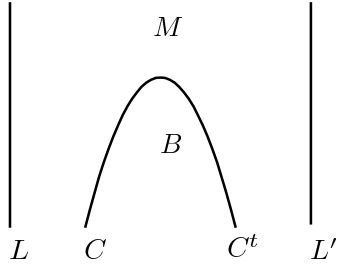}
\caption{A quilted pair of pants}\label{quiltedtriangle}
\end{figure}

Assume we have a coisotropic submanifold $C$ in a symplectic manifold $M$ which fibers over a (symplectic) manifold $B$ with sphere fibers. For example, $C$ can be a Lagrangian sphere. Denote by $\Delta_{M}$ the diagonal in $M^{-}\times M$ and let $\tau^{\pm}_{C} $ denote positive/negative fibered Dehn twist along $C$. (See Section \ref{prelimin}.)
Wehrheim and Woodward \cite{WWtriangle} (generalizing a result of Seidel \cite{seideltriangle}) prove that for any  two Lagrangian submanifolds $L,L'$ of $M$ satisfying suitable monotonicity conditions, the Lagrangian Floer chain complex $CF(L,\tau_{C} L')$   is quasi-isomorphic to the mapping  cone of the map
$$CF(L,C^{t},C,L')[\dim_{\C} B] \stackrel{\mu}\lra CF(L,\Delta_{M},L')\cong CF(L,L') $$
where the left hand side is quilted Floer homology and the map $\mu$ is given by counting \pse quilted triangles (quilted pairs of pants ) as in Figure \ref{quiltedtriangle}.
Here $C^{t}C$ and $\Delta_{M}$ can be regarded as zero and one ``resolutions'' of the fibered Dehn twist $\tau_{C}$. 
It follows from the invariance and duality properties of Floer homology that 
$CF(L,\tau^{-}_{C} L')$ is quasi-isomorphic to the cone of
$$\Bigl( CF(L,\Delta_{M},L')  [\dim_{\C} B] \stackrel{\mu^{t}}\lra CF(L,C^{t},C,L') \Bigr) [-1]$$
where $\mu^{t}$ is induced by the transpose of the \pse quilts that give $\mu$. %In Theorem \ref{myversiontri} we give a simplified description of the map $\mu$.

 Given a collection $C_{1},C_{2},\ldots, C_{N}$ of spheric coisotropic submanifolds of $M$ (where $C_{i}$ fibers over a manifold $B_{i}$) and a vector of signs $\E=(\eps_{1},\eps_{2},\ldots, \eps_{N})$, we combine the result of Wehrheim and Woodward with the above principle to obtain a spectral sequence converging to 
\bq\label{target-intro}
HF(L,\tau^{\eps_{N}}_{C_{N}}\circ\tau^{\eps_{N-1}}_{C_{N-1}}\circ\cdots \circ \tau^{\eps_{1}}_{C_{1}} (L'))
\eq
 from the $N$-dimensional hypercube obtained by 
$2^{N}$ ways of resolving the $\tau_{C_{i}}$.  
More precisely the hypercube is given as follows. If $\eps_{i}=1$ set $C_{i}^{0}=C_{i}^{t}C_{i}$ and $C_{i}^{1}=\Delta_{M}$ which are both generalized correspondences from $M$ to itself. 
 Otherwise set $C_{i}^{0}=\Delta_{M}$ and $C_{i}^{1}=C_{i}^{t}C_{i}$. For %$I=(\delta_{1},\ldots,\delta_{N})\in \{0,1\}^{N}$ 
$I=(I_{1},I_{2},\ldots, I_{N}) \in \{0,1\}^{N}$
set
\begin{equation}\label{vertexI}
CF_{I}= CF(L,C_{N}^{I_{N}}, C_{N-1}^{I_{N-1}},\ldots, C_{1}^{I_{1}}  ,L') 
\end{equation}
which is quilted Lagrangian Floer chain complex. (See section \ref{mapsofcube} for definitions.) Denote $\sigma(I)=\sum_{I_{i}=0} \dim_{\C} B_{i}$ and let $n_{-}$ denote the number of negatives among the $\epsilon_{i}$.
 As a chain group the hypercube is given by 
\begin{equation}\label{hypercube}
CF_{\oplus}=\bigoplus_{I\in\{0,1\}^{N}} CF_{I}[\sigma(I)-n_{-}] .
\end{equation} 
The maps between the adjacent vertices of the hypercube are given by counting rigid \pse quilted triangles. 
In general there are nonzero maps $\mu_{I,J}$ between nonadjacent vertices  of the cube corresponding to $I,J$. These maps are given by counting special families of \pse quilted polygons. See Section \ref{hypercubesect} for details.
%More precisely we have the following result.
%
\textit{
In this paper we work entirely with Floer homology over $\Z/2$.}

\bthm\label{mainthm}
Assuming $M$, $L,L'$ and the $C_{i}$ satisfy the Admissibility Assumption of Definition \ref{admisscond}, there is a finite cubic filtration
 on the chain complex of \eqref{target-intro}.
Therefore there is a spectral sequence  converging to \eqref{target-intro} whose first page is given by
\begin{equation}
E_{1}%=HF_{\oplus}:
=\bigoplus_{I} HF(L,C_{N}^{I_{N}}, C_{N-1}^{I_{N-1}},\ldots, C_{1}^{I_{1}}  ,L') [\sigma(I)-n_{-}]
\end{equation}
with $d_{1} $  being the sum of the maps $\mu_{I,J}$, between adjacent vertices, given by counting quilted pairs of pants.
%$=\sum_{I,J} \mu_{I,J}$
%where $J$ is an immediate successor of $I$ and $\mu_{I,J}$ is given by counting \pse quilted triangles.
\ethm

See Theorem \ref{thequasiiso} for a more precise statement and the proof. %Let $|I|=N-\sum_{i} I_{i}$.
The spectral sequence exists more generally for $L,L'$ generalized Lagrangian submanifolds of $M$ and is natural with respect to equivalence of Lagrangian correspondences (Prop. \ref{naturalfunct}).

\bcor\label{inequalities}
With the same assumptions as in Theorem \ref{mainthm}, if Lagrangian Floer homology groups are $\Z/n$ graded ($n=\infty$ for $\Z$ grading) then 
\begin{equation}
\dim HF^{i}(L,\tau^{\eps_{N}}_{C_{N}}\circ\tau^{\eps_{N-1}}_{C_{N-1}}\circ\cdots \circ \tau^{\eps_{1}}_{C_{1}} L')
\leq 
\sum_{j}\sum_{\sigma(I)-n_{-}+j=i \mod n}\dim HF^{j}(L,C_{N}^{I_{N}}, C_{N-1}^{I_{N-1}},\ldots, C_{1}^{I_{1}}  ,L').
\end{equation}
If the Floer homology groups are not graded the inequalities still hold with for ungraded homology ($n=0$).\\
\ecor

In the case where the $C_{i}$ are Lagrangian spheres (or equivalently $B_{i}$'s are points), we have
\begin{equation*}
HF(L,C_{N}^{I_{N}}, C_{N-1}^{I_{N-1}},\ldots, C_{1}^{I_{1}}  ,L')\cong 
HF(L,C_{k_{1}})\tens HF(C_{k_{1}},C_{k_{2}})\tens\cdots\tens HF(C_{k_{m-1}},C_{k_{m}})\tens HF(C_{k_{m}},L')
\end{equation*}
where $m\leq N$ and $k_{i}$'s  are so that $C_{k_{i}}^{I_{k_{i}}}\neq \Delta$.
   Also   $d_{1}$ is given by the count of \pse triangles either with one (for $\epsilon_{i}=1$) or two (for $\epsilon_{i}=-1$) outgoing ends. One feature of our spectral sequence which seems to be new, is that it involves holomorphic quilts (e.g. polygons) with multiple outgoing ends (for negative twists).
 %the composition in the Fukaya category of $M$.  
%A special case is when we have iteration of a Dehn twist, i.e. when $C_{N}=C_{N-1}=\cdots=C_{1}$.\\ %See \cite{} for a related computation. 

A typical situation where we have composition of Dehn twists, and so iterated mapping cones, is when we have a representation of the braid group on the symplectomorphism group of a manifold e.g. in symplectic Khovanov homology of Seidel and Smith \cite{SS}.
In particular we obtain a spectral sequence converging to symplectic Khovanov homology from Khovanov homology (with $\Z/2$ coefficients). (See section \ref{khsexample}.)
%We also, in section \ref{wwexample}, obtain a similar spectral sequence from Khovanov homology to the link invariant of Wehrheim and Woodward. 

\subsection{Application to \HG homology}

Fibered Dehn twists also give rise to actions of the mapping class groups of surfaces by symplectomorphisms. An important   example of this kind is given by Perutz's (and Lekili's) approach to \HG homology. Let $\Sigma_{g}$ be a surface of genus $g$.
Perutz \cite{PerutzI} assigns to each embedded circle $\gamma\subset \Sigma_{g}$ a Lagrangian correspondence $V_{\gamma}$ between the symmetric products  $Sym^{g}\Sigma_{g}$ and $Sym^{g-1}\Sigma_{g-1}$ (with appropriate symplectic forms). 
Let $H,H'$ be two handlebodies with $\partial H=\partial H'=\Sigma_{g}$. 
%of genus $g$
% given  by the alpha curves in Figure \ref{abcurves} 
%and let $H'$ be any other handlebody of genus $g$.
% $g$-tuples of circles  $\alpha_{1},\ldots,\alpha_{g}$ and $\beta_{1},\ldots,\beta_{g}$ in $\Sigma_{g}$.
 %which constitute a symplectic basis for $H_{1}\Sigma_{g}$ (as in Figure \ref{abcurves}). 
Denote by $T_{H}$ and $T_{H'}$ their corresponding Heegaard tori.
\bthm\label{HGss}
Let $M$ be a closed oriented 3-manifold obtained by gluing $H$  to $H'$ by a homeomorphism $\phi$ of the surface $\Sigma_{g}$.
Let $\phi$ be given as a composition of Dehn twists along curves in $\Sigma$,
$\phi=\tau_{\gamma_{N}}^{\eps_{N}}\circ \cdots \circ \tau_{\gamma_{1}}^{\epsilon_{1}}$. Choose a basepoint $z\in\Sigma$ away from the curves $\alpha_{i}, \beta_{i}$ and $\gamma_{i}$. 
Then there is a spectral sequence which converges to $\HFh(M)$ and whose $E_{1}$ term is given by
\begin{equation}\label{hge1}
\bigoplus_{I} HF(T_{H'}, V_{\gamma_{N}}^{I_{N}},\ldots,V_{\gamma_{1}}^{I_{1}}, T_{H}).
\end{equation}

\ethm

In the above theorem the symplectic manifolds involved are $Sym^{g}\Sigma_{g}\bs (\{z\}\times Sym^{g-1}\Sigma_{g})$ and $Sym^{g-1}\Sigma_{g-1}\bs (\{z\}\times Sym^{g-2}\Sigma_{g-1})$ with the same K\"ahler forms used in Perutz construction.
 One also Hamiltonian isotopes the submanifolds $V_{\gamma_{i}}$ to become balanced. The differential $d_{1}$ is again given by counting \pse quilted pairs of pants.
Note that since $c_{1}(Sym^{g}\Sigma_{g})$ is nonzero (even mod $n$), the above homology groups are ungraded.

Theorem \ref{HGss} is interesting even in the case $N=1$ where it is a reincarnation of the knot surgery exact triangle \cite{BranchedDouble} in terms of the mapping class group. 
 Each summand in \eqref{hge1} is the $\HFh$ of a 3-manifold given as follows. Let $Y_{\gamma_{i}}$ denote the cobordism between $\Sigma_{g}$ and $\Sigma_{g-1}$ in which $\gamma_{i}$ is excised out and let $Y_{\gamma_{i}}^{I_{i}}$ be defined with the same rule as for $V_{\gamma_{i}}^{I_{i}}$ (eg. if $\epsilon_{i}=1$, and $I_{i}=0$ then $Y_{\gamma_{i}}^{I_{i}}=Y_{\gamma_{i}} Y_{\gamma_{i}}^{t} $ ).
Then it follows from the work in progress of Lekili and Perutz \cite{lekili-p} that the  $I$'th summand of \eqref{hge1} is the $\HFh$ of the 3-manifold given by the sequence of cobordisms
$H, Y_{\gamma_{N}}^{I_{N}},\cdots, Y_{\gamma_{1}}^{I_{1}}, H' $.\\
 %  by a sequence of $N$ cobordisms $Y_{i}$ connecting $H$ to $H'$ such that $Y_{i}$ is either the trivial cobordism or the one in which $\gamma_{i}$ is collapsed to a point and then restored back.

%We do not attempt at giving a combinatorial description of \eqref{hge1} in this paper  however it is feasible.
We compute the $E_{1}$ page in the case that $M$ is the double cover of $S^{3}$ branched over a link $L$ (denoted $\Sigma(K)$) and show that the $E_{1}$ page is Khovanov's hypercube for $L$ (with coefficients in $\Z/2$). More precisely we have the following.

\bthm\label{osthm}
%With the same assumptions as 
In Theorem \ref{HGss} let $\phi$ be hyperelliptic, so  its mapping cylinder is the double cover of $S^{2}\times [0,1]$ branched over a braid $b\in Br_{2g}$. Also assume $H'=H$ be a handlebody of genus $g$.  
%be the branched double cover of a set of $g$ unlinked arcs $\{c_{i}\}_{i=1}^{g}$. 
Then the  page $(E_{1},d_{1})$ of the spectral sequence of Thm.~\ref{HGss} is naturally isomorphic to the Khovanov hypercube (over $\Z/2$) for the plat closure $K$ of the mirror of $b$ and its $E_{\infty}$ page is $\HFh(\Sigma(K)\,\#\, S^{1}\times S^{2})$.
% by the $c_{i}$ where $\Sigma(K)$ is the branched double cover of $K$.  
\ethm

%The theorem holds for other types of closure (by arbitrary crossingless matchings) as well.
 The proof relies on work in progress by Lekili and Perutz \cite{lekili-p}. (See \ref{commrel} below.) 
The spectral sequence of Thm.~\ref{osthm} is very closely related to the spectral sequence of Ozsvath and Szabo \cite{BranchedDouble}. %However it has two minor differences. Firstly the $E_{2}$ page gives Khovanov homology (over $\Z/2$) and not reduced Khovanov homology. 
However the (higher) maps in our hypercube are not exactly given by (higher) cobordism maps for \HG homology.
Lipshitz, Ozsvath and Thurston \cite{lipshitz:branchedI} have been working on obtaining an explicit computation of the spectral sequence of Ozsvath and Szabo using bordered \HG homology. Our approach can provide an alternative, more geometric method which we hope to get back to in future. 
% In fact the higher differentials in  the spectral sequence of Thm.~\ref{osthm} reduce to polygon maps on the genus one surface. 
 %Note that for general $M$, the $E_{2}$ term of the spectral sequence of Theorem \ref{HGss} is not a topological invariant. See Remark \ref{watsonresult}. 

\begin{remark}\label{watsonresult}
It may be tempting to call the $E_{2}$ page of the spectral sequence of Theorem \ref{HGss} the ``Khovanov homology of the 3-manifold $M$``.
However it has been shown by Watson \cite{watson:khbranched} that Khovanov homology of a link $K$ does not give an invariant of the branched double cover of $K$. This, together with Thm. \ref{osthm}, implies that  this $E_{2}$ page % of the spectral sequence of Theorem \ref{HGss} 
  is not a 3-manifold invariant. Nonetheless in a forthcoming work we study this $E_{2}$ term in more detail and give a combinatorial description for it.
\end{remark}

\begin{question}
Which 3-manifolds have presentations for which the spectral sequence of Theorem \ref{HGss} collapses at $E_{2}$ page?
\end{question}

It follows from the work of Ozsvath and Szabo \cite[Prop 3.3]{BranchedDouble} together with Thm.~\ref{osthm} that branched double covers of quasi-alternating links have this property. It is possible that the manifolds in question are an appropriate generalization of L-spaces. (L-spaces are  rational homology 3-spheres $M$ for which $\dim_{\Z/2}\HFh(M)=\# H_{1}(M,\Z)$). 

 Corollary \ref{inequalities} now gives upper bounds on the rank of $\HFh(M)$ from each presentation of $M$ by gluing of handlebodies.
One possible application of these inequalities 
%of in the case of \HG homology 
is in deciding whether a given 3-manifold can be obtained from a given composition of (classical) Dehn twists.
We end this introduction by noting that obtaining such a spectral sequence for a topological invariant, defined using Lagrangian Floer homology, can be regarded as the first step toward obtaining a combinatorial description of the invariant. 
%\red{In fact John Baldwin and Adam Levine  obtinined a new combinatorial model for knot Floer homology using a method whih can be compared to the general method of this paper. }

\begin{remark}
The same arguments that are used to prove the Theorem \ref{mainthm} can be adapted to prove a hypercube in the derived Fukaya category $D\F^{\#}(M,M)$. This means that, under the same assumptions as in Theorem \ref{mainthm}, the Lagrangian correspondence $\on{graph} (\tau^{\eps_{N-1}}_{C_{N-1}}\circ\cdots \circ \tau^{\eps_{1}}_{C_{1}})$ is isomorphic, in $D\F^{\#}(M,M)$, to an element of the form $(\sum_{I} (C_{N}^{I_{N}},\cdots, C_{1}^{I_{1}}),\, D)$ where $D$ is given by the count of quilts similar to the ones that give the differential on the hypercube.
See Proposition \ref{opencube}.
\end{remark}

%\subsection{Further questions for \HG homology}
%Theorem \ref{HGss} raises a few questions.

%\begin{question}
%How is the spectral sequence of Theorem \ref{HGss} related to the work of Manolescu and Ozsvath \cite{manolescu:hyperbox}?
%\end{question}
%
%There are (at least) two methods for constructing all oriented  3-manifolds: one by Dehn surgery along links and the other by gluing handlebodies using surface homeomorphisms. For Theorem \ref{HGss} we use the latter approach while Manolescu and Ozsvath use the former. However the construction of  \cite{manolescu:hyperbox} is more general in that it investigates how the \HG homology of a general 3-manifold changes under a Dehn  surgery.

\subsection*{Organization}
In section \ref{prelimin} we recall basic facts about spheric coisotropic submanifolds.
In section \ref{examples} we recall some examples of coisotropic submanifolds of importance to low dimensional topology.
Section \ref{hypercubesect} is the technical part of the paper in which we construct the hypercube under different admissibility conditions on the manifolds involved.
The exact triangle for fibered Dehn twists is reviewed in section \ref{wwtrisection}.
In section \ref{spectrsect} we prove Theorem \ref{mainthm}.
Theorems \ref{HGss} and \ref{osthm} are proved in sections \ref{hgproof} and \ref{osproof} respectively.

\subsection*{Acknowledgements.}  
I  thank Paolo Ghiggini for his lectures on \HG homology and for a helpful discussion. I would also like to thank Robert Lipshitz for pointing out the relevance of Watson's result and Tim Perutz for  helpful correspondence.
A part of this work was done when the author was visiting the Mathematical Sciences Research Institute as a research member. The author was also partly supported by QGM (Center for Quantum Geometry of Moduli Spaces).% funded by the Danish National Research Foundation.

\section{Preliminaries}\label{prelimin}
The main objects of study in this paper are fibered Dehn twists along spherically fibered\\ coisotropic submanifolds. 
%
%A theorem of Perutz and shows that such fibered Dehn twists are Hamiltonian isotopic to monodromy maps of this fibration.
Such submanifolds are closely related to 
 symplectic Morse-Bott fibrations (also called Lefschetz-Bott fibrations) which  are generalizations of symplectic Lefschetz fibrations. 
%
%We will be working in a restricted version of the symplectic category in which morphisms are given by spherically fibered coisotropics.  
%
Let $(M,\omega)$ be a symplectic manifold and $C\subset M$ coisotropic (which means that for any $x\in C$ and any $v\in T_{x}M$, if $\omega(v,w)=0$ for all $w\in T_{x}C$ then $v$ is tangent to $C$). Because of the closedness of $\omega$, the distribution $\ker \omega|_{C}$ is integrable and the resulting foliation is called the \emph{null foliation} of $C$. 
\bdf
A coisotropic submanifold $C$ of a symplectic manifold $(M,\omega)$ is  \emph{fibered} if there is a manifold $\bar{M}$ and a fibration $\pi:C\to \bar{M}$ which is constant on the leaves of the null foliation. 
It is \emph{spherically fibered} (or just \emph{spheric}) if the fibers of $\pi$ are spheres $S^{k}$ and moreover the structure group of the bundle $\pi$ can be reduced to $SO(k+1)$ where $k$ is the codimension of $C$. 
\edf
The manifold $\bar{M}$ inherits a symplectic from $\omega$ given by
\begin{equation}
\bar{\omega}(\pi_{*}v,\pi_{*}w)=\omega(v,w).
\end{equation}
Let $\iota:C\to M$ be the inclusion then $(\iota,\pi):C\to M^{-}\times \bar{M}$ is a Lagrangian embedding and so we can regard $C$ as a Lagrangian correspondence from $M$ to $\bar{M}$ (or equivalently $C^{t}=(\pi,\iota)C$ a correspondence from $\bar{M}$ to $M$).

%Fibered coisotropic submanifolds are closely related to (singular) symplectic fibrations and in particular spherically fibered ones are related to Lefschetz-Bott fibrations. 
 %With notation as in the last paragraph 
$SO(k+1)$ acts on $\C^{k+1}$ with a moment map $\eta: \C^{k+1}\to \so (k+1)^{*}$ whose regular fiber is $S^{k}$.
The $SO(k+1)$-bundle associated to $\pi$ yields an associated $\C^{k+1}$-bundle $P$ over $\bar{M}$. 
Let $\alpha \in \Omega^{1}(P,\so(k+1))$ be a connection one form on $P$ and let $\pi_{1},\pi_{2}$ be the projections from $P$ to $\bar{M}$ and $\C^{k+1}$ respectively.
One has a  closed  2-form on $P$ given by
\begin{equation}
\omega'= \pi_{1}^{*}\bar{\omega}+\pi_{2}^{*}\omega_{\C^{k+1}}+ d\langle\eta,\alpha\rangle
\end{equation}
where $\langle,\rangle:\so(k+1)^{*}\tens \so(k+1)\to \R$ is the paring.
This form  is nondegenerate in a \nbhd $P_{\epsilon}$ of the zero section \cite{guillemin:inphysics}.
\renewcommand{\b}{\mathbf}

Consider  $T^{*}S^{k}$ as a subset of $\C^{k+1}$ given by the pairs $(\mathbf{x},\mathbf{y})$ such that 
$|\mathbf{x}|=1, \langle \b{x},\b{y}\rangle=0 $. Let $\psi\in C^{\infty}(\R,\R)$ be monotone increasing on $[0,1]$, send this interval to $[\pi,2\pi]$ and be equal to  $\pi$ or $2\pi$ outside this interval. The \emph{generalized Dehn twist} $\tau$ along the zero section is  defined by
\begin{equation}
\tau\left(\begin{array}{c}
\b{x}\\ \b{y}
\end{array}
 \right)
=\left( \begin{array}{cc}
\cos(\psi(|\b{y}|))\cdot I & |\b{y}|^{-1}\sin(\psi(|\b{y}|))\cdot I\\
-|\b{y}|\sin(\psi(|\b{y}|))\cdot  I     &   \cos(\psi(|\b{y}|))\cdot I
 \end{array} \right)
\left(\begin{array}{c}
\b{x}\\ \b{y}
\end{array}
 \right)
\end{equation}
where $I$ is the $(k+1)\times (k+1)$ identity matrix. If $\sigma_{t}$ is the Hamiltonian flow of the length function $|\b{y}|$ then $\tau(\b{x},\b{y})=\sigma_{\psi(|\b{y}|)} (\b{x},\b{y})$.
It is easy to see that $\tau$ extends smoothly to the zero section and is the antipodal map there.
One can see by direct computation that $\tau$  is a symplectomorphism for the canonical symplectic structure on $T^{*}S^{k}$. 
% \red{as it is the Hamiltonian diffeomorphism induced by the function}. 

Now  there is a symplectomorphism $\tilde\tau$ of $(P_{\epsilon},\omega')$ which is a generalized Dehn twist along the $S^{k}$ in each $\C^{k+1}$ fiber.
More precisely $\tilde\tau$ is the time $2\pi$ map of the Hamiltonian flow of $\psi\circ \eta$.
 %More precisely since $P$ is given as an associated bundle, any diffeomorphism of $\C^{k+1}$ induces a diffeomorphism on $P$.
%
 %It can be given explicitly as the Hamiltonian diffeomorphism associated to a function on $P$ (which we leave to the reader).
By the coisotropic \nbhd theorem \cite[Thm 39.2]{guillemin:inphysics} a \nbhd of $C$ in $M$ is symplectomorphic to a \nbhd of the zero section in $P_{\epsilon}$. This way $\tilde\tau$ induces a symplectomorphism $\tau_{C}$ of $M$ called the \emph{fibered Dehn twist} along $C$. The Hamiltonian isotopy class of $\tau_{C}$ is not changed by the action of Hamiltonian isotopies on $C$. However it is not known if this class is independent of the ``framing'' i.e. the choice of a local symplectomorphism into $P_{\epsilon}$.

%Moreover the monodromy of this fibration is  Hamiltonian isotopic to the fibered Dehn twsit along $C$.

%\newcommand{\Sympr}{$\mbf{Symp}^{r}$}
%\red{The morphisms can also be graphs of dehn twists. The definition may not be needed.}
%\bdf
%The restricted symplectic category \Sympr has symplectic manifolds as objects. Morphisms between two objects $M$ and $N$ are given by equivalence classes of pairs of sequences $(\{L_{i} \}_{i=1}^{n},\{ M_{i}\}_{i=0}^{n})$ such that $M_{0}=M, M_{n}=N$ and $L_{i}\subset M^{-}_{i-1}\times M_{i}$ is  a spherically fibered coisotropic. 
%The equivalence relation between morphisms is given as follows. 
%Morphisms are composed by concatenation. 
%\edf

% Morphisms in \Sympr  are called generalized Lagrangian correspondences. 
%A (restricted) \emph{generalized Lagrangian submanifold} of a symplectic manifold $M$ is a morphism in \Sympr from $M$  to the point  given by $(\{L_{i} \}_{i=1}^{n},\{ M_{i}\}_{i=0}^{n})$ \st $\dim M_{i}>\dim M_{i+1}$ for each $i$.

% \red{ is any such coisotropic the vanishing cycle of a LB fibration?} 

\section{Some examples of spheric coisotropic submanifolds}\label{examples}

\bx[following Seidel and Smith \cite{SS} ]\label{ssex}
Let $\mathcal{S}_m\subset \sl_{2m}(\C)$ be the set of matrices   of the form  \\ \bq\label{matrixS}\left(
                                         \begin{array}{cccccc}
                                           y_1 & I &\space & & \space& \\
                                           y_2 & & I &  & & \\
                                           \vdots&  &  & \ddots & \\
                                             y_{n-1}& & & & I\\
                                             y_n &  &  & & 0\\
                                         \end{array}
                                       \right)\eq\\
where $I$ is the $2\times 2$ identity matrix, $y_1\in \mathfrak{sl}_2$ and $y_i\in \mathfrak{gl}_2$ for $i>1$.
$\mathcal{S}_{m}$ is a transverse slice to the adjoint orbit of a nilpotent matrix of Jordan type $(m,m)$ in $\sl_{2m}$ \cite[Lemma 23]{SS}. 
Let $\Sigma_{2m}$ denote the symmetric group on $2m$ letters and
consider the map $\chi:\S_{m}\to \C^{2m}/\Sigma_{2m}$ which sends each matrix to its spectrum. Set 
\begin{equation}
\Y_{m,\nu}=\chi^{-1}(\nu).
\end{equation}
If $\nu$ has no repetitions then it is a regular value of $\chi$  and so $\Y_{m,\nu}$ is a K\"ahler (in fact Stein) manifold with  vanishing first Chern class. 
 %(We can regard $\C$ as a subset of the Riemann sphere and so non compactness of $\C$ is not a problem.)
Denote by $Conf_{2m}$ the set of regular values of $\chi$, i.e. $2m$-tuples without multiplicities.
Let $\delta:[0,1)\to Conf_{2m}$ be any curve such that $\nu':=\lim_{t\to 1} \delta(t)$ has an element $\mu$ of multiplicity two and no more repetitions among its members. Let $\bar{\nu}\in Conf_{2m-2}$ be the result of deleting $\mu$ from $\nu'$.
It can be  easily shown that the set of critical points of $\chi$ in $\Y_{m,\nu'}$ is symplectomorphic to $\Y_{m-1,\bar{\nu}}$ so we can regard $\Y_{m,\nu'}$ as a submanifold of $\Y_{m,\nu'}$. \\

Let $D\subset$  be a small disk containing the image of $\delta$ and such that $D\cap Conf_{2m}=D\bs \{\delta(1)\}$. Seidel and Smith prove \cite[Lemma 27]{SS} that up to a K\"ahler isomorphism the restriction of $\chi$ to $D$ is of the form %$\pi$ where %$\pi:\C^{3}\to \C$ is given by 
$\pi(x,a,b,c)=a^{2}+b^{2}+c^{2}$.
Let $U\subset \Y_{m-1,\bar{\nu}}$ be compact.
Let $L_{\delta}=L_{\delta}(U)\subset \Y_{m,\mu}$ be the set of points of $\Y_{m,\delta(t)}$ (for $t$ close to $1$) which converge to an element of $U\subset \Y_{m-1,\bar{\nu}}\subset \Y_{m,\nu'}$  under  the gradient flow of $Re \pi$. 
% is defined for all $t\in[0,1)$ and converge to elements of the critical point 
$L_{\delta}$ is a relative vanishing cycle for the fibration $\chi|_{D}$.
  Morse lemma  implies that after possibly replacing $\nu$ with $\delta(t)$ for $t$ close to $1$, the map $L_{\gamma}\to \Y_{m-1,\bar{\nu}}$ which sends a point to its limit under the gradient flow is smooth and is a $S^{2}$-bundle over its target. So, $L_{\delta}$ is a spheric coisotropic submanifold of $\Y_{m,\nu}$ which fibers over $U\subset \Y_{m-1,\bar{\nu}}$. Moreover if we choose a larger compact subset $U'$, the restriction of the resulting  $L_{\delta}(U')$ to $U$ is Hamiltonian isotopic to $L_{\delta}(U)$.

%\grn{Maybe we can use Ciprian's alternative description also }\\ 

This construction has been generalized by Manolescu \cite{ManolLinkInvs} to the case of $\sl_{n\cdot m}$. The analogs of the above coisotropic submanifolds are $\C P^{n}$ bundles in that case.\\

\ex

\bx[following Perutz  \cite{PerutzI}]\label{perutzex}
Let $S$ be a Riemann surface of genus $g$ and $\gamma\subset S$ an embedded circle. Let $\bar{S}$ be the result of surgering  $\gamma$ out and attaching  two discs to $S\bs \gamma$.
 Let $\eta\in H^{2}(Sym^{g}S)$ be Poincare dual to $\{pt\}\times Sym^{g-1} S$ and $\theta$ be such that $\theta-g\cdot\eta$ is Poincare dual to 
$\sum a_{i}\times b_{i} \times Sym^{g-2}S$ where 
%$g$ is the genus of $S$ and
 $\{a_{i},b_{i}\}$ is a symplectic basis for $H_{1} (S)$.
Let $P_{S}\in H^{2} Sym^{g}S$ be a cohomology class of the form $t\eta+s\theta$ where $s,t$ are constants.
As in \cite{perutz:hamilslides} there are symplectic forms in this class which agrees with the push-forward of the product symplectic form (on $S\times S\times \cdots \times S$) away from the big diagonal.  
 We similarly have a cohomology class $P_{\bar{S}}$ for $Sym^{g-1} \bar{S}$ (with the same values of $s$ and $t$). 

\bthm[Perutz \cite{PerutzI}, Theorem A]\label{perutzthm}
If $\omega$ and $\bar{\omega}$ are symplectic forms in cohomology classes $P_{S}$ and $P_{\bar{S}}$ respectively then 
there is a coisotropic submanifold $\iota: V_{\gamma} \to (Sym^{g}S,\omega)$ which fibers over $(Sym^{g-1}\bar{S},\bar{\omega})$ with circle fibers.
%\times (Sym^{k-1}\bar{S},\bar{\omega}) $. Moreover the projection $\pi:V_{\gamma}\to  is a circle bundle.

\ethm
As in Example \ref{ssex},
$V_{\gamma}$ is the vanishing cycle for a Lefschetz-Bott fibration  over the disk.  The generic fiber of this fibration, which we denote by $p$, is $Sym^{g}S$ and  its set of critical points can be identified with $Sym^{g-1} S$. More specifically one starts with a Lefschetz fibration $\pi$ over the disk in which $S$ becomes nodal along $\gamma$ and then one considers the fibration $p$ where the $p^{-1}(z)$ for each point $z$  is the 
Hilbert scheme of points on $\pi^{-1}(z)$. 
This Hilbert scheme for nonsingular curves (i.e. $\pi^{-1}(z)$ for $z\neq 0$) is the same as the symmetric product. 
See \cite{PerutzI} for details.

\ex

\bx[following Wehrheim and Woodward \cite{WWfieldb}]\label{wwex}
%This example is of type II.
% For a surface $S\in \M_{g,2k+1}$, let $M(S)$ denote the moduli space of flat $SU(2)$ connections on $S$ whose holonomy around each puncture of $S$ lies in the conjugacy class of ``traceless'' matrices. $M(S)$ \red{is smooth and has a symplectic structure given by }

Let $\Sigma_{g,n}$ denote a topological surface of genus g with n punctures and let $\M_{g,n}$ denote the moduli space of flat SU(2) connections on $\Sig_{g,n}$ whose holonomy around each puncture has trace zero.  If $n$ is odd then this moduli space is smooth \cite[Prop. 3.3.1]{WWfieldb}.
$\M_{g,n}$ also has a symplectic structure which goes back to Atiyah and Bott. See \cite{alekseev:groupvalued} for a more modern approach.
 
Let $\gamma\subset \Sig_{g,n}$ be an embedded circle  which does not bound a (punctured) disk and whose complement is connected.  
Then one can consider the three dimensional cobordism $Y_\gamma$ between $\Sig_{g,n}$ and $\Sig_{g-1,n}$ in which $\gamma$ is pinched to a point and then excised out. 
Consider $C_\gamma\subset \M_{g,n} \times \M_{g-1,n}$ consisting of pairs of connections which extend to the whole of $Y_\gamma$. Note that projection on the first factor embeds $C_{\gamma}$ in $\M_{g,n}$.
 For two nearby punctures $z_{1}, z_{2}$ one can also consider the cobordism $Y_{z_{1},z_{2}}$ between $\Sigma_{g,n}$ and $\Sigma_{g,n-2}$ in which $z_{1}$ and $z_{2}$ merge.  It gives rise to a subset $C_{z_{1},z_{2}}\subset \M_{g,n}$.

\bthm[Wehrheim-Woodward \cite{WWfieldb}, Prop. 3.4.2]
 $C_{\gamma}$ is a smooth coisotropic submanifold of $\M_{g,n}$;  it fibers over $\M_{g-1,n}$ with $S^3$ fibers. The set $C_{z_{1},z_{2}}$ is also a smooth coisotropic submanifold of $\M_{g,n}$ which fibers over $\M_{g,n-1}$ with $S^{2}$ fibers.
%3) It is a Lagrangian submanifold of $\M_g,n^- \times M_g-1,n$ where minus sign means that you multiply the symplectic form by -1
\ethm

%A Lagrangian submanifold of $M^-\times N$ is called a Lagrangian correspondence between $M$ and $N$. 
So an elementary cobordism between surfaces gives a Lagrangian correspondence between the corresponding moduli spaces of flat connections.
%%%%%
Because $C_\gamma\subset \M_{g,n}$ is coisotropic, one can consider the fibered Dehn twist $\tau_{C_\gamma}$ along it. 
A result of Wehrheim and Woodward \cite{WWtriangle}  (extending results of Callahan and Seidel) shows that fibered twist along $C_\gamma$ agrees up to Hamiltonian isotopy, to the diffeomorphism induced on $\M_{g,n}$ by the classical Dehn twist along $\gamma$.
%It follows that the assignment $\gamma \to\tau_{C_\gamma}$   gives a homomorphism from the mapping class group of $\Sig_{g,n}$ to the $\pi_{0}$ of the  symplectomorphism group of $\M_{g,n}$.

\ex

\section{The hypercube}\label{hypercubesect}
In this section we construct the maps between the vertices of the hypercube of \eqref{hypercube} and show that they make \eqref{hypercube} into a chain complex. These maps are given by the count of \pse quilts.
The conditions that guaranty that the maps in the hypercube are well-defined and give rise to a differential turn out to be the same as the conditions that guaranty the Fukaya category $\F(M)$ of the symplectic manifold $M$ is well-defined.
%
%Let $M$ be a compact symplectic manifold which has a well-defined Fukaya category $\F(M)$. 
%Description of Fuk: Floer theoretic info about lags
Bubbling of \pse discs with boundary on a single Lagrangian    or \pse polygons with fixed Maslov index but with unbounded energy prevent the Fukaya category from being well-defined. The latter problem is prevented by imposing the monotonicity (or exactness) condition on the Lagrangians while for the first problem there are a number of different remedies. Below we recall a few well-known situations in which $\F(M,\omega)$ is well-defined.

\begin{enumerate}
\item \label{wwfuk}
 $M$ is monotone and one includes only the Lagrangians $L$ which are monotone, the image of $\iota:\pi_{1}(L)\to \pi_{1}(M)$ is torsion and the minimal Maslov number of $L$ is greater than two. % \cite[Lemma 4.1.3]{QuiltedFloer}.
\item \label{seidelfuk}
 $\omega$ is exact with convex (contact type) boundary and one includes only its compact exact Lagrangian submanifolds (i.e. if $\omega=d\sigma$ then $\sigma|_{L}$ is exact) \cite{Seidelbook}.
\item \label{balancedfuk}
 $M$ is monotone with a prequantum line bundle $K$ with a connection form $\alpha$ and one includes only  balanced (also called Bohr-Sommerfeld monotone) Lagrangians $L$ (i.e. $K|_{L}$ has a section $s$ for which $s^{*} \alpha$ is exact) with the additional condition that $\pi_{2}(M,L)=0$. % \cite[Lemma 4.1.5]{QuiltedFloer}.
%
%\item \label{absmithfuk}
%$M$ is a Calabi-Yau 4-manifold (so, $2c_{1}(M)=0$) and one includes only monotone Lagrangians $L$ with vanishing Maslov class \cite[Lemma 3.3]{AbSmithI}.
%
\item \label{myfuk}
$M$ is a (noncompact) Stein manifold of finite type and one includes only exact Lagrangian submanifolds which are invariant under the Liouville flow outside a compact subset 
\cite[Section 4.3]{RR1}.
\end{enumerate}

\bdf[\textbf{Admissibility Condition}]\label{admisscond}
% In the rest of this paper
    A symplectic manifold $M$ and a collection $(L,C_{1},\ldots, C_{k},L')$ where $L,L'\subset M$ are Lagrangian and $C_{i}\subset M$ are spheric coisotropics fibering over manifolds $B_{i}$
are said to be \emph{admissible} if they satisfy one of the above conditions ($C_{i}$ as a Lagrangian submanifold of   $M^{-}\times B_{i}$). For \eqref{wwfuk} and \eqref{balancedfuk}  the $C_{i}$ and $L,L'$ are assumed to have the same monotonicity constant.
\edf
%In each of these cases the generalized Fukaya category $\F^{\#}(M)$ is well-defined if one includes only generalized correspondences
%$M\stackrel{L_{1}}\to M_{1}\stackrel{L_{2}}\to M_{2}\to\cdots \to pt$
%\st each $M_{i}$ and each $L_{i}\subset M^{-}_{i-1}\times M_{i}$ satisfy the  condition (with the same monotonicity constant). For \eqref{myfuk}  the components must be proper as well.

\bl\label{nobubble}
If $(M,\omega)$ and $L\subset M$ are admissible then for any almost complex structure $J$ compatible with $\omega$, any $J$-holomorphic map from the disk to $M$ that sends the boundary of the disk to $L$ is constant.
\el

This follows from exactness in \eqref{seidelfuk} and \eqref{myfuk} and by the vanishing of relative $\pi_{2}$ in \eqref{balancedfuk}. %For \eqref{absmithfuk} it is shown in  \cite[Lemma 3.3]{AbSmithI}.
For \eqref{wwfuk} this is shown e.g. in \cite[Thm. 1.2]{OhFloer}.

\subsection{A family of quilts}\label{familyofquilts}
\renewcommand{\Q}{\mathcal{Q}}
%\red{We don't need to worry about one point of a cap getting close to a 'strip point' because we can fix 3 points.}

In this section we introduce a family of quilts  and in the next  we use this family to define the maps on the hypercube.
%Quilts that we use %to define the hypercube are a special class of those studied by Wehrheim and Woodward 
%
For more on quilts see \cite{WWquilts}.

Let $N>0$ be a fixed integer and $\E=(\eps_{1},\ldots,\eps_{N})\in \{-1,1\}^{N}$. Put the lexicographic ordering on $\{0,1\}^{N}$. 
Let  $ I,J \in \{0,1\}^{N}$ be \st $I\leq J$.
%we define a family of quilts $\mathcal{Q}_{I,J}$ as follows.
%
%
Let $R_{I,J}$ be the set of isomorphism classes of  pairs $(S_{I},S_{J})$ where $S_{I}$ is the unit circle with a  fixed point $z_{+}:=\sqrt{-1}$, called the \emph{outgoing point}, and a finite set of marked points distinct from $z_{+}$ whose number is determined by $I$ and $J$ as follows. Give the circle minus $z_{+}$ the counterclockwise orientation.
%If $I_{i}=J_{i}=0$, a single point on the circle minus $z_{+}$ corresponds to $i$ and is called a \emph{solitary point}.
 If $I_{i}<J_{i}$ and $\eps_{i}=+1$, a couple of adjacent points is associated to $i$ and they are called a \emph{pair}.
This is subject to the following conditions.  
\begin{itemize}
\item If $i<j$ then  the points corresponding to $i$ are before those corresponding to $j$ in the counterclockwise orientation.
\item If $I_{i}<J_{i}$ and $I_{j}< J_{j}$ with $j=i+1$ 
%and between points corresponding to $i$ and $j$ there are only solitary points (if any)
then the left marked  point corresponding to $j$ and the right marked point corresponding to $i$ 
%and all the solitary points in between
 are identified.
\end{itemize}

The marked circle $S_{J}$ has a similar structure. The distinguished point on its boundary is called the \emph{incoming point} and denoted by $z_{-}=-\sqrt{-1}$. If $\eps_{i}=-1$ and $I_{i}<J_{i}$, a couple of points on the circle minus $z_{-}$ are assigned to $i$
% and if $I_{i}=J_{i}=1$ a solitary point
 with the same identification and ordering requirement (now clockwise).  $R_{I,J}$ is the quotient of the the set of such pairs of marked circles with the action of the Mobius transformations of the disc.

\newcommand{\Rb}{\overline{R}}

Since $R_{I,J}$ is a moduli space of points on two copies of the circle, it has the Deligne-Mumford-Stasheff compactification  $\Rb_{I,J}$ by stable tree-discs as in Section 9f of \cite{Seidelbook}. Codimension one components of the boundary of $R_{I,J}$ correspond to two marked points (not necessarily adjacent) in $S_{I}$ or $S_{J}$ (but not both at once)   converging together.\\

\begin{figure}
\includegraphics[width=0.8\textwidth]{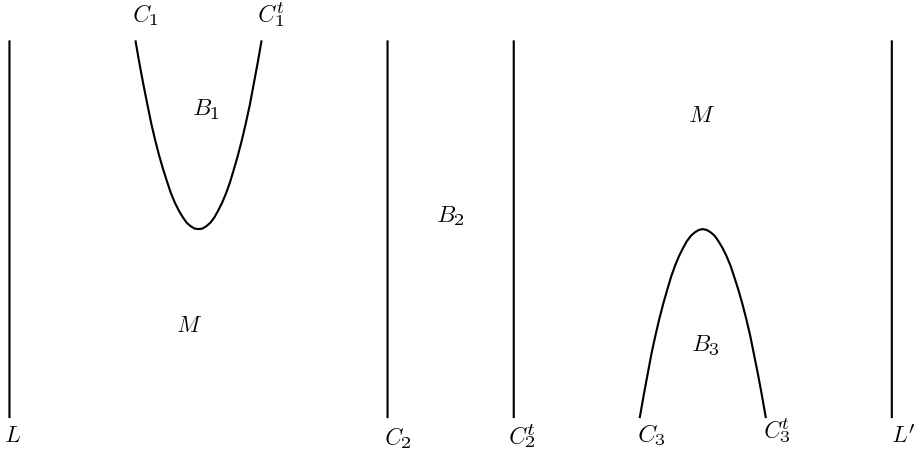}
\caption{The quilt $Q_{I,J}$ for $\E=(-1,1,1)$ and $I=(0,0,0), J=(1,0,1)$  with labelling.}\label{quiltexample}
\end{figure}

\begin{figure}
\includegraphics[scale=0.7]{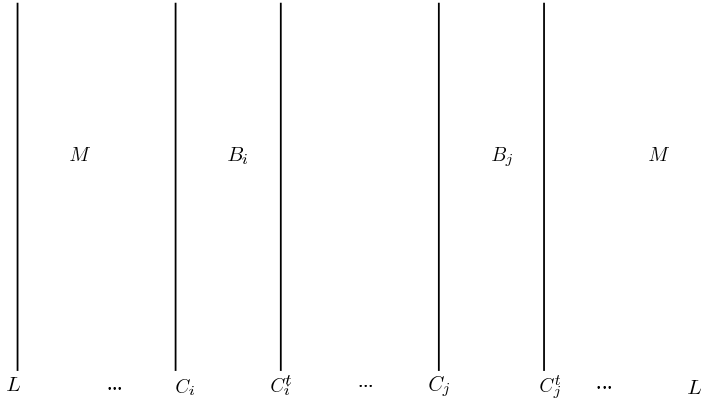}
\caption{A quilted strip}\label{quiltedstrip}
\end{figure}

%Each quilt $Q \in \Q_{I,J}$ 
To $R_{I,J}$ corresponds a quilt   $Q_{I,J}$ which can be described as follows.
%\red{clarify what a strip is.}
$Q_{I,J}$ can be thought of as lying in $[0,2N+1]\times \R\subset \C$ and satisfies the following; which give all the seams of $Q$.
%\begin{itemize} 
%\item 
If $\eps_{i}=+1$ then $Q\cap ([2i-1,2i]\times \R)$ contains a semi-infinite cap if $I_{i}<J_{i}$ 
and is an infinite strip\footnote{i.e. $ \{2i-1\}\times \R$ and $ \{2i\}\times \R$ are seams in $Q$.} 
 if $I_{i}=J_{i}=0 $. %  no strips associated to $\Delta_{M}$
If $\eps_{i}=-1$ then $Q\cap ([2i-1,2i]\times \R)$ contains a semi-infinite cup if $I_{i}<J_{i}$ and is an infinite strip if $I_{i}=J_{i}=1$.
%
%\item The $y$ coordinate of the tip of each cap is negative and is less than or equal to that of caps closer to the boundary of $[0,2N]\times \R$. Similarly the $y$ coordinate of the tip of each cup is positive and is greater than or equal to that of caps closer to the boundary of $[0,2N]\times \R$. 
%\end{itemize}
%
 Figure \ref{quiltexample} shows an example of such a quilt. For such a $Q$ let $Q^{0}$ denote $[0,2N+1]\times \R$ minus all the cups, caps and strips mentioned above. %It is diffeomorphic to a disjoint union of strips.
%$\Q_{I,J}$ is the quotient of this family by the free $\R$ action acting by vertical translation. 
Note that for any $I$, $\Q_{I,I}$ consists of a single quilted strip as in Figure \ref{quiltedstrip}.

%For each $I\leq J$ the family $\Q_{I,J}$ can be compactified as follows.

One assigns quilts to the boundary strata of $R_{I,J}$ as follows.
 Let $r$ be a point in the codimension one stratum of $\partial R_{I,J}$.  %corresponding to the converging of the marked points corresponding to $I'$.
 There is always a subinterval\footnote{In the obvious sense.} %$I'\subset I$ (or $J'\subset J$) 
$A\subset \{1,2,\ldots,N\}$
 \st $r$ corresponds to the converging of the marked points corresponding to $A$ in $S_{I}$ or $S_{J}$. 
Consider the case of $S_{I}$.
 Let $I'\in \{0,1\}^{N}$ be obtained from $I$ as follows.
 We set $I'_{i}=I_{i}$ if $i\notin A$.
For $i\in A$, if  $I_{i}=J_{i}$ we set $I'_{i}=I_{i}$ and if $I_{i}<J_{i}$ we set $I'_{i}=J_{i}$. 
It is evident that $I< I' < J$.
 To $r$ we assign the pair of quilts $(Q_{I,I'}, Q_{I',J})$.  You can see that gluing $Q_{I,I'}$ to the bottom of $Q_{I',J}$ gives back $Q_{I,J}$. The construction for the case that $r$ corresponds to $S_{J}$ is similar.
The construction for the strata of  codimension $k$ is similar and gives a $(k+1)$-tuple of quilts (organized by a tree).\\
%\red{tree-quilts} similar to the tree-discs of  Section 9e of \cite{Seidelbook}.

% Let $\{Q_{n}\}\subset \Q_{I,J}$ be a sequence of quilts \st the $y$ coordinates of the tips of some of the caps in them go to $+\infty$ and/or the $y$ coordinates of some of the cups go to $-\infty$, the limiting quilt in $\bar{\Q}_{I,J}$ has infinite strips in place of those caps and cups. 
 %Note also that the codimension one part of this boundary is given by two cups (possibly the same) and all the cups in between them go to negative infinity or two caps (possibly the same) and all the caps in between them go to positive infinity (but not both caps and cups at the same time) \red{Why codim 1?}.

Since the quilts we consider  have the same combinatorial type as polygons, the same method as that of  Lemma 9.3 in \cite{Seidelbook} can be used to deduce the existence of a consistent choice of strip like ends for quilts associated to the strata of $\Rb_{I,J}$. This means that there exists a compact family $\Q_{I,J}$ of quilts and a projection
\bq 
 \pi : \Q_{I,J} \to \overline{R}_{I,J}
\eq
 such that $\pi^{-1}(r)$ for any $r\in \overline{R}_{I,J}$ is a quilt with the combinatorial type described above. (So, for example for $r\in R_{I,J}$, $\pi^{-1}(r)$ has the same type as $Q_{I,J}$.) %Moreover
Concretely this involves choosing, for each $r\in \Rb_{I,J}$, the  heights of the caps in cups
% and strips
 in the quilt  associated to $r$ in such a way that when $r$ converges to a  point $r'$, the limiting heights agree with the ones given to the $\pi^{-1}(r')$ by gluing.
One can also use Theorem 1.3 in \cite{MWW} for this purpose.

Note that since the points of $\Rb_{I,J}$ are marked circles up to Mobius transformations, the points of  $\Q_{I,J}$ are  quilts up to vertical translation (i.e. Mobius transformations fixing positive and negative infinity).
So for $N=1$ %because of quotienting by the $\R$ action, 
the family of quilts is compact. Indeed for $I=(0)$ and $J=(1)$, $\Q_{I,J}$ consists of a single quilted triangle (pair of pants) as in  Figure \ref{quiltedtriangle}.
% This is true also for the subset quilts in $\Q_{I,J}$ (for any $N$) that have only one cup or cup. 

\subsection{The maps of the hypercube}\label{mapsofcube}

As before let $M$ be a symplectic manifold, $L,L'\subset M$ Lagrangian submanifolds and $C_{1},\ldots,C_{N}$ spheric coisotropic submanifolds of $M$ \st $C_{i}$ fibers over a symplectic manifold $B_{i}$.  
We label the seams of the quilts  $Q_{I,J}$ as follows.   The boundary seams $\{0\}\times \R$ and $\{2N+1\}\times \R$ are labeled by $L$ and $L'$ respectively.
 For $i\in \{1,2,\ldots,N\}$ if $\{2i-1\}\times \R$ and $\{2i\}\times \R$ (or a subinterval of them) are  part of  seams in $Q$ then they are labeled by $C_{N-i+1}$.

Recall \cite{QuiltedFloer} (also \cite[ section 3.1]{perutz:gysin}) that a \emph{generalized Lagrangian correspondence} (from point to itself)
is a sequence
\bq\label{alagcor} 
\mathbf{L}=\bigl(pt \stackrel{L_{1}}\to M_{1} \stackrel{L_{2}}\to M_{2}\to\cdots \to M_{k}\stackrel{L_{k+1}}\to pt \bigr)
\eq
 of symplectic manifolds $(M_{i},\omega_{i})$ and Lagrangian correspondences between them $L_{i}\subset M_{i-1}^{-}\times M_{i}$.
By adding a diagonal correspondence $\Delta_{M_{i}}\subset M_{i}^{-}\times M_{i}$ if necessary we can assume that the number $k$ is odd. Set
$\mathcal L_{0}=L_{2}\times L_{4}\times\cdots\times L_{k+1}$ and
$\mathcal \L_{1}=L_{1}\times L_{3}\times \cdots\times L_{k}$ which are Lagrangian submanifolds of $\prod_{i} M_{i}$ (with the symplectic form $\sum (-1)^{i} \omega_{i}$ ).
We say that the Lagrangian correspondence \eqref{alagcor} \emph{intersects transversely} if $\mathcal L_{0}$ and $\mathcal L_{1}$ intersect transversely. If this is the case and the $L_{i}$ are compact then the intersection $\mathcal L_{0}\cap \mathcal L_{1}$ is a finite set and its elements are the \emph{generators} of the quilted Floer chain group $CF(\mathbf{L})$. These generators can be described as $k$-tuples
$(x_{1},\ldots,x_{k})\subset M_{1}\times M_{2}\times\cdots\times M_{k}$
\st $(x_{i},x_{i+1})\in L_{i+1}$.\\

Let $I,J\in \{0,1\}^{N}$ be \st $I\leq J$. In order to define the maps between the vertices \eqref{vertexI} of the hypercube, corresponding to $I$ and $J$
we first need to choose Hamiltonian isotopies that make each one of the  generalized correspondences $(L,C_{N}^{I_{N}}, C_{N-1}^{I_{N-1}},\ldots, C_{1}^{I_{1}}  ,L')$ and  $(L,C_{N}^{J_{N}}, C_{N-1}^{J_{N-1}},\ldots, C_{1}^{J_{1}}  ,L')$ intersect transversely. 
We also need to choose for each $Q_{I,J}\in \Q_{I,J}$ a family of compatible almost complex structures on $M$ and the $B_{i}$, parameterized by the points\footnote{More precisely parameterized by the $x$ coordinate of the points outside a compact subset.} of each patch of $Q_{I,J}$
 which make the linearization of the equation of \pse maps surjective (and so yield smooth moduli spaces). In addition such choices for  $Q_{I,J}$, $Q_{J,K}$ and $Q_{I,K}$ with $I\leq J\leq K$ must be consistent under gluing.
This is done in an inductive way from highest codimension strata of $\Q_{I,J}$ (which consist of quilted strips and/or quilted pairs of pants) up to its interior. See Theorem 6.24 in \cite{MWW}. Such a set of Hamiltonian isotopies and almost complex structures  is called a \emph{regular set of perturbation data}.
 In the sequel we assume that such Hamiltonian isotopies and almost complex structures are chosen and are used to define the set of generators and the equation for \pse curves.

For each $I,J$ with $I\leq J$ denote by $\M_{I,J}$ the moduli space of pairs $(Q,\ul{u})$
where $Q\in \Q_{I,J}$ and $\ul{u}=(u_{0},\ldots , u_{N})$ 
is \st
\begin{itemize}
\item for $j>0$, if $Q\cap [2j-1,2j]\times \R$ is a cap, a cup or a strip, $u_{j}$ is a \pse map from that cap/cup/strip to $B_{j}$ (the base of the fibration of the coisotropic $C_{j}$) and is a space holder otherwise.
\item $u_{0}: Q^{0}\to M$ is \pse.
\item Whenever there is a seam between two components of $Q$ and $u_{j},u_{k}$ are their corresponding \pse maps then the pair $(u_{j},u_{k})$ sends the seam to the Lagrangian correspondence labeling the seam.  \\
\end{itemize}

%By the general theory of \pse quilts \cite[]{}, 
Outside a compact set each quilt is a quilted strip ( as in Figure \ref{quiltedstrip}) which is the same thing as a strip in the product manifold with boundary on the product Lagrangians $\L_{0},\L_{1}$. Therefore \pse quilts have the exponential convergence property which means that 
each such $\ul{u}$ maps $-\infty$ to a generator 
%$\ul{x}$
 of $CF_{I}$ and  $+\infty$ to a generator
% $\ul{y}$
 of $CF_{J}$. So, 
$\M_{I,J}$ is a disjoint union of $M_{I,J}(\ul{x},\ul{y})$ for $\ul{x}, \ul{y}$ generators of $CF_{I}$ and $CF_{J}$ respectively. 
$\M_{I,J}(x,y)$ gets a topology from the topology of $R_{I,J}$ and the fact that for $i>0$ (resp. $i=0$), $u_{i}$  lies in the $W^{1,2}$ space of maps from the strip/cap/cup (resp. $Q^{0}$) to $B_{i}$ (resp. $M$). The regular choice of almost complex structures on $M$ and the $B_{i}$ implies that $\M_{I,J}(x,y)$ is a disjoint union of  smooth manifolds
\footnote{More precisely this means that $M_{I,J}(x,y)$ is locally a smooth submanifold of the Cartesian product of $R_{I,J}$ with the space of $W^{1,2}$ maps.} of possibly different dimensions.
Let $\M_{I,J}(\ul{x},\ul{y})_{d}$ denote the $d$ dimensional part of  $\M_{I,J}(\ul{x},\ul{y})$.
If $\ul{u}=(u_{1},\ldots, u_{N})$ then by definition the energy of $\ul{u}$ is the sum of the energies of the $u_{i}$.
\bl\label{energybounded}
%Let $\Q$ be a family of quilts each one one of them has one incoming and one outgoing end labeled by two fixed Lagrangian correspondences $\ul{L}_{0}$ and $\ul{L}_{1}$. If $ul{L}_{i}$ satisfy the Admissibility Condition then the curves in the zero dimensional part of the moduli space have bounded energy.
For each $d$ there is a constant $k_{d}$ \st the energy of each $\ul{u}$, for $(Q,\ul{u})$ element  of $\M_{I,J}(x,y)_{d}$, is less than $k_{d}$.
\el
 \bp
For exact Lagrangians this follows easily from Stokes theorem applied to the pullback of the appropriate symplectic forms by the $u_{i}$.
For other types we note that as in the proof of Theorem 3.9 in \cite{WWquilts}, near each interior seam the couple of maps $(u_{i},u_{j})$ (where either  $j=0$ or $j=i+1$) is equivalent to a map $u'_{i}$ from the strip into the product $M^{-}\times B_{i}$  which sends one boundary component to $C_{i}$. This is the case because the seams are assumed to be real analytic (and we are negating the symplectic form on $B_{i}$). Note that for the quilts that we use this neighborhoods are dense in each component of the quilt. 
This way, two  $\ul{u}$, $\ul{v}$ give, for each $i$, a map from the cylinder $\S^{1}\times \R$ to $M^{-}\times B_{i}$ whose area is the difference of the energies of $u'_{i}$ and $v'_{i}$. Therefore the lemma will follow if we show that the energy of any \pse map  $w: S^{1}\times [0,1]\to M^{-}\times B_{i}$, which sends $S^{1}\times \{1\}$ to $C_{i}$, is proportional to the Maslov index of the curve $w|_{S^{1}\times \{1\}}$ i.e. $E(w)=k\cdot \mu(w|_{S^{1}\times \{1\}})$ where $k$ is the same for all the $C_{i}$ and $L,L'$.
 This is because this Maslov index
%  obtained from $u'_{i}$ and $v'_{i}$
 has to be less than a number which does not  depend on the curves otherwise $\ul{u}$ and $\ul{v}$ would not be in the $d$ dimensional part of the moduli space. 

%the sum of the Maslov indices gives the local dimension of the   moduli space $\M_{I,J}(x,y)$.
%and therefore the curves in $\M_{I,J}(x,y)_{d}$,

If $L,L'$ and the $C_{i}$ are admissible of type \eqref{wwfuk},  the assumption on the fundamental groups implies that each such $w$ is homotopic to a map of the disk and then one can use the monotonicity of $C_{i}, L$ or $L'$.
If the Lagrangian and coisotropic submanifolds   admissible of type \eqref{balancedfuk}, this is shown in \cite[Lemma 4.1.5]{QuiltedFloer}
\ep

\bl\label{mycompact}
If $M,L,L'$ and the $B_{i}$ satisfy the Admissibility Condition \eqref{myfuk} then the elements of $\M_{I,J}(x,y)$ lie in a compact subset of $M$.
\el
\bp
It is shown in \cite[Lemma 3.3.2]{RR2}, using an argument of Oh \cite{Ohnoncompact}, that holomorphic curves with boundary on two exact Lagrangian submanifolds of a Stein manifold satisfying admissibility condition of type \eqref{myfuk}, lie in a compact set determined (only) by the intersection points of the Lagrangians. One generalizes this to \pse quilts using the same ``folding'' argument used in Lemma \ref{energybounded}.
\ep

\bl \label{M0finite}
If $M,L,L'$ and the $B_{i}$ satisfy the Admissibility Condition  then $\M_{I,J}(x,y)_{0}$ is a finite set.
\el
\bp
Let $(Q_{n},\ul{u}_{n})$ be a sequence of elements of $\M_{I,J}(x,y)_{0}$. If $Q_{n}$ converges to an element of the boundary of $\Q_{I,J}$, which for simplicity we assume to be in the codimension one part  and so of the form $(Q_{I,I'},Q_{I',J})$ for $I<I'<J$,  then because of the boundedness of energy Lemma \ref{energybounded}, $u_{n}$ will converge, up to reparametrizaion, to a pair
 % in $C^{\infty}$ topology, 
  $(u_{I,I'}, u_{I',J})$ where $u_{I,I'}\in \M_{I,I'}$ and $u_{I',J}\in \M_{I',J} $. 
 So by the gluing theorem for \pse quilts ( Proved in \cite[Theorem 1]{Maugluing} for \pse quilted polygons) there is a one parameter family of \pse quilts and this, together with the fact that $\M_{I,J}(x,y)$ is locally a manifold, contradicts the assumption that $(Q_{n},u_{n})$ lie being in the zero dimensional part of $\M_{I,J}(x,y)$.
The same argument prevents the limit from being a  ``broken quilt'' i.e. of the form $(u_{I,I},u_{I,J})$ or $(u_{I,J},u_{J,J})$.
%Knowing this, the rest of the proof is the same as that of Theorem 3.9.b in \cite{WWquilts}. 
%To elaborate, Gromov compactness generalizes to \pse quilts \cite[]{} and implies that \red{Why bounded energy?} the limit of the sequence $u_{n}$ is a \pse quilt with possible bubbling or a broken quilts. Bubbling is ruled out by the Admissibility Condition.
% and broken quilts are impossible due to zero dimensionality assumption. 
The only other possible limit is a bubbling of a boundary disk which is ruled out by Lemma \ref{nobubble}.
Therefore $\M_{I,J}(x,y)_{0}$ is compact.
\ep

%  By \cite[]{}, for each $Q\in \Q$,  the zero dimensional part is finite and 
Therefore one can define a linear map $\mu_{I,J}:CF_{I}\to CF_{J}$ given on the generators of $CF_{I}$ by

\begin{equation}
\mu_{I,J}(\ul{x})=\sum_{\ul{y}} \# \M_{I,J}(\ul{x},\ul{y})_{0} \cdot \ul{y}
\end{equation}
where the sum is over the generators of $CF_{J}$.
Note that $\mu_{I,I}$ is the (quilted) Floer differential on $CF_{I}$. Also if the manifolds $B_{i}$ are  points and the $\eps_{i}$ are positive then the maps $\mu_{I,J}$ are higher composition maps in the Fukaya category $\F(M)$ of $M$. 

Define a linear map $D: CF_{\oplus}\to CF_{\oplus}$ by
\begin{equation}
D=\sum_{I\leq J} \mu_{I,J}.
\end{equation}
Let $G_{I}$ denote the set of generators of $CF_{I}$.
\bl\label{d2=0}
If $M$, $L,L'$ and the $C_{i}$ satisfy the Admissibility Condition then $D^{2}=0$.
\el

\bp
This is a special case of the Master Equation for quilt families \cite[Theorem 1.5]{MWW}.
We note that if $\ul{x}\in G_{I}$ then 
\begin{equation} 
D^{2}\ul{x}=\sum_{J\geq I} \sum_{K\geq J} \sum_{\ul{y}\in G_{J}}\sum_{\ul{z}\in G_{K}} \#\M_{I,J}(\ul{x},\ul{y})_{0}\times \M_{J,K}(\ul{y},\ul{z})_{0} \cdot \ul{z}.
\end{equation}
 Therefore if we show that $\partial \M_{I,K}(\ul{x},\ul{z})_{1} $  can be identified with
\bq\label{boundaryofM1}
 \bigcup_{I\leq J\leq K}\bigcup_{\ul{y}\in G_{J}} \M_{I,J}(\ul{x},\ul{y})_{0}\times M_{J,K}(\ul{y},\ul{z})_{0}
\eq
 the result follows. The proof is standard and similar to that of Lemma \ref{M0finite}. If $(Q_{n},u_{n})$ is a sequence in $\M_{I,K}(\ul{x},\ul{y})_{1}$ then, by the gluing argument, the limit of $\{Q_{n}\}$ can not lie in the part of $\partial R_{I,J}$ of codimension two or higher because then the $u_{n}$ would not be in the one dimensional part of $\M_{I,K}(\ul{x},\ul{y})$. So, the limit of $Q_{n}$, if not in the interior of $R_{I,J}$, is a pair $(Q_{I,J},Q_{J,K})$ with $I\leq J\leq K$. Therefore $u_{n}$ converges, up to reparametrization, to a pair $(u_{I,J},u_{J,K})$ where $u_{I,J}\in \M_{I,J}$ and $ u_{J,K}\in \M_{J,K}$ have to be in the zero dimensional part.
 It follows from exponential convergence that $u_{I,J}$ (resp. $u_{J,K}$) sends $+\infty$ (resp. $-\infty$) to an element of $G_{J}$ which have to be the same.
%
%or corresponding to a limiting point of the quilt family $\Q_{I,K}$ which gives an element of $\M_{0}(x,y)\times M_{0}(y,z)$ \red{WHY?}
%
Bubbling is again ruled out by \ref{nobubble}. 
Therefore $\partial \M_{I,K}(\ul{x},\ul{y})_{1}$ is included in  \eqref{boundaryofM1}.
%$\cup_{I\leq J\leq K}\cup_{y\in C_{J}} \M_{I,J}(\ul{x},\ul{y})_{0}\times M_{J,K}(\ul{y},\ul{z})_{0}$. 
The reverse inclusion follows from the gluing argument for \pse quilts \cite{Maugluing}.
\ep

\section{The Wehrheim-Woodward exact triangle}\label{wwtrisection}

The main tool that we use in this paper to prove Theorem \ref{thequasiiso} is the following result of Wehrheim and Woodward. See also \cite{perutz:gysin} and \cite{biran:gysin} for related results.
 Let $M$ be a symplectic manifold and $L,L'\subset M$ be Lagrangian. 

\bthm[Wehrheim-Woodward\cite{WWtriangle}, Theorem 5.2.9]\label{WWtrithm}
Let $C\subset M$ be a spheric coisotropic submanifold 
%of codimension at least two
 fibering over a manifold $B$. If the triple $(C,L,L')$ is monotone and each one of $C,L,L'$ has Maslov index at least 3 then $CF(L,\tau_{C},L')$ is quasi-isomorphic to the cone of
\begin{equation}
CF(L,C,C^{t},L')[\dim_{\C} B] \stackrel{\mu}\lra CF(L,L')
\end{equation}
where $\mu$ is given by the count of \pse quilted triangles as in Figure \ref{quiltedtriangle} i.e. $\mu=\mu_{(0),(1)}$.
\ethm 
%
%\begin{remark}
%Even though the theorem is proved with the assumption of the
%codimension of $C$ being at least two, as in \cite[17.15]{Seidelbook}
%(for the case of Lagrangian circles) the theorem holds for codimension
%equal to one over $\Z/2$.
%\end{remark}

\begin{figure}
\includegraphics{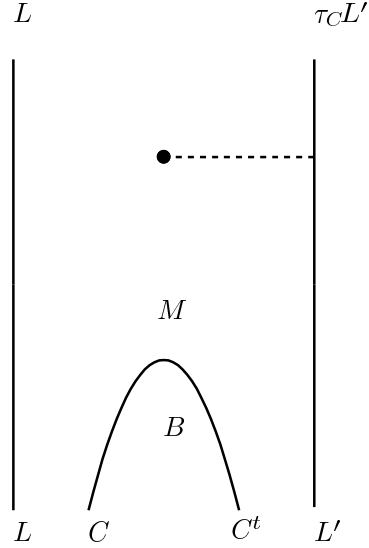}
\caption{The quilted fibration used to construct the mas $C_{1} \to
  C_{\infty}$. One considers a (unique up to isotopy) Lefschetz-Bott fibration over each patch of the  quilt
  whose fiber over any point is either $M$ or $B$ as indicated. The
  dot represents the unique critical value of the fibration.  The
  fibration is such that its monodromy around the critical point
  is Hamiltonian isotopic to $\tau_{C}$.  A \emph{section} of such quilted fibration  consists of a
  pair of
\pse sections (of the fibrations over each patch) whose value on each
seam lie in the given Lagrangian submanifold. }\label{qiso+0}
\end{figure}

\begin{figure}[ht]
\includegraphics{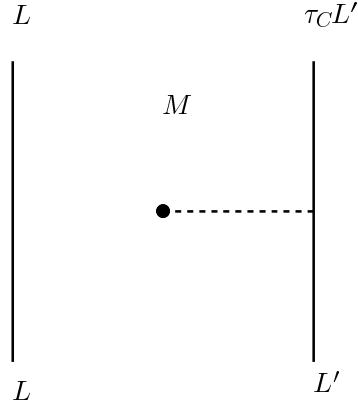}
\caption{The quilted fibration used to construct the map $C_{0} \to
  C_{\infty}$.}\label{qiso+1}
\end{figure}

The theorem more generally holds for $L,L'$ generalized Lagrangian submanifolds of $M$ (i.e. generalized correspondences between $M$ and a point). See \cite[Prop. 5.16]{RR2}.
%As shown in \cite[Corollary 5.17]{RR2}, for negative Dehn twists
We now recall the construction of the quasi-isomorphism. % and also describe a similar quasi-isomorphism for negative Dehn twists.
For simplicity denote $C_{0}=CF(L,C,C^{t},L')$, $C_{1}=CF(L,L')$ and
$C_{\infty}= CF(L,\tau_{C},L')$. It is easy to see that, in general, a chain map
from $Cone(f)$ to $C_{\infty}$ consists of a chain map $k:C_{1}\to
C_{\infty}$ and a homotopy $h:C_{0}\to C_{\infty}$ between $k\circ \mu$
and zero. Moreover $(h,k)$ is a quasi-isomorphism if and only if its
mapping cone is acyclic.
The maps $h,k$ are given by counting \pse sections of ``quilted
 Lefschetz-Bott fibrations'' in Figures \ref{qiso+0} and \ref{qiso+1}.
The map $k$ is given by the count of the rigid sections of
the quilted fibration of Figure \ref{qiso+0}.

Let $Q_{t}$ for $t\in [0,1]$ be the one parameter family of quilted fibrations \st $Q_{0}$ is the fibration in
Figure \ref{qiso+1} and as $t\to 1$, the critical
value  approaches the cap and a circle is pinched off the cap.
 The map $h$ is given by the
count of zero dimensional part of the moduli space $\{(Q_{t},s)| t\in [0,1]\}$
where $s$ is a section of the fibration $Q_{t}$. By definition (and because bubbling is ruled out) the map $h$ is
 a homotopy between $k\circ \mu$ and  the map associated
to  $Q_{1}$ which we denote by $h_{1}$.
% i.e. the quilt in Figure \ref{qiso+1} in which the distance between the critical value and the cap is zero.
 Wehrheim and Woodward \cite[Lemma 5.2.1]{WWtriangle} show that if
the codimension of $C$ is $\geq 2$ then
$h_{1}$  is zero and hence $(h,k)$ is a chain map from
$Cone(\mu)$ to $CF(L,\tau_{C},L')$.

In fact, with coefficients in $\Z/2$, this holds for the codimension one case as well. To see this we recall the proof of the vanishing of $h_{1}$.
Let $E_{C}$ denote the Lefschetz-Bott fibration over the disk whose vanishing cycle is $C$. (Such a fibration always exists; see \cite[2.4.1]{PerutzI}.) 
The quilted fibration $Q_{1}$ that gives $h_{1}$ is the result of gluing $E_{C}$ to another quilted Lefschetz-Bott fibration.
It follows from gluing theorem (along a seam) for \pse quilts that if the moduli space of \pse sections of $E_{C}$ has sufficiently high dimension near any point then the zero dimensional part of the moduli space of sections of $Q_{1}$ is empty.  This is indeed what Wehrheim and Woodward do i.e. (as in \cite{seideltriangle}) they explicitly produce a high dimensional moduli of sections of $E_{C}$ containing a given section. 
 %
%The obstruction for $h_{0}$  being zero is given by the moduli space of sections of the Lefschetz-Bott fibration obtained by gluing the fibration over the disk whose vanishing cycle is $C$ to an arbitrary quilted Lefschetz-Bott fibration. 
 For the codimension one
case, as in \cite[17.15]{Seidelbook} (for Lefschetz fibrations), the zero dimensional part of the moduli space of sections of such a fibration consists of two points and hence $h_{0}$ is zero  over $\Z/2$. The same  result should hold over $\Z$ with a careful choice of signs as in the aforementioned reference however we do not use this.

%Now we describe the chain maps for negative fibered Dehn twists. Here we have $C_{0}=CF(L,L')$, $C_{1}=CF(L,C^{t},C,L)$ and $C_{\infty}=CF(L,\tau^{-1}_{C}L')$ and $f^{t}:C_{0}\to C_{1}$ is given by the count of rigid \pse maps of the transpose of the quilt on Figure \ref{quiltedtriangle}.
%We want a chain map $(h,k)$ from the
%cone of $f^{t}$ to $C_{\infty}$. 
%The map $k$ is now given by the count of sections of the quilted
%fibration in Figure \ref{qiso+0} but with $\tau^{-1}_{C}L'$ instead of
%$\tau_{C} L'$. 
%In turn $h$ is given by the count of the moduli space %$(Q'_{t},s)$
%where $Q'_{t}$ is a quilted fibration as in Figure \ref{qiso-} and $t$ denotes the distance between the critical point and the disk. $s$ is a rigid \pse section of $Q'_{t}$. Again $h$ is a homotopy between $k\circ f^{t}$ and the map associated to $Q'_{0}$ which is zero by the same argument as for the case of positive twists. Therefore $(h,k)$ is a chain map. 

%\begin{figure}
%\includegraphics{qiso-}
%\caption{The quilted fibration used to construct the map %$CF(L,L')\to
%  CF(L,\tau_{c}^{-1} L')$.}\label{qiso-}
%\end{figure}

To prove the acyclicity of $Cone(h,k)$ as in \cite{seideltriangle} Wehrheim and Woodward use Floer homology over the Novikov ring (even though the Lagrangians are monotone). Such chain complexes are filtered by area. 
Let $D$ denote the differential on $Cone(h,k)$. One of the main ingredients in the proof is the 
decomposition of $D$ into two parts $D=D_{0}+D_{1}$ whose degrees have a gap in between i.e.  the degree of $D_{0}$ is in $[0,\eps)$ and the degree of $D_{1}$ is in $[2\eps,\infty)$ for some $\epsilon>0$.
The other ingredient is using a geometric argument to prove the acyclicity of $D_{0}$  for fibered  Dehn twists with small support. The acyclicity of $Cone(h,k)$ follows from these two by homological algebraic arguments.

%\red{More precisely such a \pse quilted triangle is a pair $(u_{1},u_{2})$ \st}

As a consequence one has an exact triangle
\[
\xymatrix{&\ar[dl]_{[1]} HF(L,\tau_{C}L') &\\
          HF(L,C^{t},C,L')[\dim_{\C} B] \ar[rr]_{\mu}& & HF(L,L')\ar[ul]_{k}
}
\]
where the map $\stackrel{[1]}\to$ is the connecting homomorphism.

\begin{remark}\label{negtri}
One has a corresponding exact triangle for negative Dehn twists. We have $CF^{i}(L,\tau^{-1}_{C}L')\simeq CF^{i}(\tau_{C}L,L')\simeq \Hom(CF^{l-i}(L',\tau_{C}L),\Z/2) $ where $l=\dim L$. This induces an isomorphism on Floer homology. Noting that the map induced by the dual $Q^{t}$ of a quilt $Q$ is the dual of the map induced by $Q$, Theorem \ref{WWtrithm} yields a quasi-isomorphism from $CF(L,\tau^{-1}_{C}L')$ to 
$Cone(\mu^{t})[-1]$ where $\mu^{t}:CF(L,L')\to CF(L,C^{t},C,L')$ is given by counting the upside-down version of the quilt in Figure \ref{quiltedtriangle}.
%the cone of $\mu^{t}:CF(L,L')\to CF(L,C^{t},C,L') $. 
This results in an exact triangle as follows.
\[
\xymatrix{& HF(L,\tau^{-1}_{C}L')\ar[dr]^{k^{t}} &\\
          HF(L,C^{t},C,L')[\dim_{\C} B-1] \ar[ur]^{[1]}& & HF(L,L')[-1]\ar[ll]^{\mu^{t}}
}
\]
\end{remark}

\bcor
If $M$ and $L,L',C$ are admissible then the conclusion of Theorem \ref{WWtrithm} holds.
\ecor
\bp
The assumption of Maslov index at least 3 in Theorem \ref{WWtrithm} is to rule
out disc bubbling. In our case bubbling is ruled out by Lemma \ref{nobubble}.
Monotonicity of $(C,L,L')$ for \eqref{wwfuk} follows from \cite[Lemma 4.1.3]{QuiltedFloer},
for \eqref{balancedfuk} from \cite[Lemma 4.1.5]{QuiltedFloer}
and for \eqref{seidelfuk} and \eqref{myfuk} from exactness.
Finally for \eqref{myfuk} we need to show that all the \pse curves involved lie in a compact subset of $M$. This was shown in \cite[proposition 5.14]{RR2}. (See Lemma \ref{mycompact}.)
\ep

Here we make a simple observation about  the map $\mu$  used in the exact triangle which is useful in computations. Let $\pi:C\to B$ be the projection. Also let $\alpha,\beta,\gamma$ be the components of the boundary of a disk with three punctures.

\bpr\label{myversiontri}
%Let $J_{z}$ be a family of almost complex structures on $M$, parametrized by the points of a triangle, which is regular for the triple $(L,L',C)$.
%Then the moduli space of \pse quilts used in defining the map $\mu$ in Theorem \ref{WWtrithm} is cobordant to the (zero dimensional part of the) moduli space of $J_{z}$-holomorphic triangles in $M$ with boundary conditions given by the Lagrangians $L,L'$ and the coisotropic $C$. 
For a regular family of almost complex structures, the map $\mu$ is given by (the zero dimensional part of) the moduli space of \pse maps from the disk into $M$ which send $\alpha,\beta,\gamma$ respectively to $L,L',C$
and
 for which $\pi\circ u|_{\gamma}$ can be completed to a \pse disk in $B$.
\epr
\bp
%Pick two families $J, J'$ of almost complex structures on $M$ and $B$ respectively which give transversal moduli space of \pse representations of the labeled quilt $Q$ n Figure \ref{quiltedtriangle} transversal.
Let $u$ and $v$ be %$J-$  and $J'-$holomorphic  
\pse maps from the two components of $Q$ into $M$ and $B$ respectively. %which satisfy the boundary conditions
By the abuse of notation let $\gamma$ denote the seam between the two components of this quilt.
Since $C$ is embedded in $M$, we have $v|_{\gamma}=\pi\circ u|_{\gamma}$ where $\pi: C\to B$ is the projection.
%$u$ determines the restriction of $v$ to the seam and
 By unique continuation for \pse maps \cite{Ahn}, $v$ is uniquely determined by $v|_{\gamma}$ and therefore by $u$. 
%Such a \pse cap $v$ does exist because $\pi\circ u|_{\gamma}$ is analytic.
\ep

\section{The spectral sequence}\label{spectrsect}
In this section we show that the chain complex of \eqref{target-intro} is quasi-isomorphic to the hypercube $CF_{\oplus}$ of resolutions of the twists constructed in the last section.
We do this by an inductive argument involving the exact triangle of Wehrheim and Woodward.
 We use the abbreviation 
\bq
CF:=CF(L,\tau^{\eps_{N}}_{C_{N}}\circ\tau^{\eps_{N-1}}_{C_{N-1}}\circ\cdots \circ \tau^{\eps_{1}}_{C_{1}} (L')).
\eq
Theorem \ref{mainthm} follows from the following.

\bthm\label{thequasiiso}
%There is a chain map $\phi: CF_{\oplus} \to CF$ which induces an isomorphism on cohomology.
Assume $L,L'$ and the $C_{i}$ satisfy the Admissibility Condition.
Then $CF_{\oplus}$ and $CF$ are isomorphic in the derived category $D^{b} ({\Z/2}\on{-mod})$. In other words there is a chain complex $CF^{+}$ and chain maps $f: CF^{+}\to CF_{\oplus}$, $g:CF^{+}\to CF$ which induce isomorphisms on homology.
Moreover if the Floer  homology groups are graded, this quasi-isomorphism preserves the grading.
\ethm

%
%\red{$\Z/2$ grading for oriented ones}

Before proceeding to the proof, we recall the generalities about grading on Lagrangian Floer homology \cite{gradedlag}. Let $(M,\omega,J)$ be a symplectic manifold with a compatible almost complex structure and assume there is an $n>0$  \st  $2c_{1}(M)$ is zero in $H^{2}(M,\Z/n)$. This implies that there is a (non-unique) line bundle $\eta$ and an isomorphism $r:\eta^{\tens n} \to \Lambda^{max} (TM,J)^{\tens 2}$.
Let $\on{Lag(M)}$ denote the fiber bundle over $M$ whose fiber over any point $x$ is the Lagrangian Grassmannian of $T_{x} M$.
The isomorphism classes of such pairs $(\eta,r)$ is in one to one correspondence with fiber bundles $\L\to M$ whose fiber over each point $x$ is an $n$ fold cover of the Lagrangian Grassmannian of $T_{x} M$ (corresponding to the Maslov class)
 together with a map  $\L\to \on{Lag}(M)$ which is an $n$ fold covering. For each Lagrangian submanifold $L\subset M$ there is a canonical section $s_{L}$ of $\on{Lag}(M)$. A grading on $L$ is a lift of this section to a section of $\L$. A choice of grading for two Lagrangians $L,L'$ induces an absolute $\Z/n$ grading on $HF(L,L')$. \\

%Let $C=Cone(f:C_{0}\to C_{1})$ and $D=Cone(g:D_{0}\to D_{1})$ be chain complexes. If $F:C\to D$ is a chain map then $F$ decomposes into a matrix $[F_{ij}]_{i,j=0,1}$ where $F_{ij}:C_{i}\to D_{j}$.

%\bl\label{homollem}
%$F_{01}$ is a homotopy between $g\circ F_{00}-F_{11}\circ f$ and zero. Moreover if $F'_{0,1}$ is any other such homotopy and $F'$ is obtained 
%\el
%Proof is a simple comutation.

The quasi-isomorphisms $f,g$ are given explicitly in terms of  maps induced by quilted fibrations.
As the first step we note that $CF$ is quasi-isomorphic to 
\begin{equation}\label{decomposedtwist}
CF(L,\on{graph} \tau_{C_{N}}^{\epsilon_{N}},\ldots, \on{graph}\tau_{C_{1}}^{\epsilon_{1}},L').
\end{equation}
The chain complex $CF^{+}$ is the hypercube associated to resolving only the positive twists, i.e. ones for which $\epsilon_{i}=+1$. 
Let $n_{+}$ denote the number of positive twists in \eqref{decomposedtwist} and for $I\in \{0,1\}^{n_{+}}$ let $C^{+}_{I}$ be the Floer chain complex for the correspondence obtained from $(L,\on{graph} \tau_{C_{N}}^{\epsilon_{N}},\ldots, \on{graph}\tau_{C_{1}}^{\epsilon_{1}},L')$ in which the $i$th positive fibered twist is replaced with its $I_{i}$-resolution.
Set
\begin{equation}
CF^{+}=\bigoplus_{I\in \{0,1\}^{n_{+}}} CF^{+}_{I}.
\end{equation}
One can define maps $\mu_{I,J}: CF^{+}_{I}\to CF^{+}_{J}$ for $I\leq J$ as before and the sum of all such $\mu_{I,J}$ makes $CF^{+}$ into a chain complex.

\begin{figure}
\includegraphics{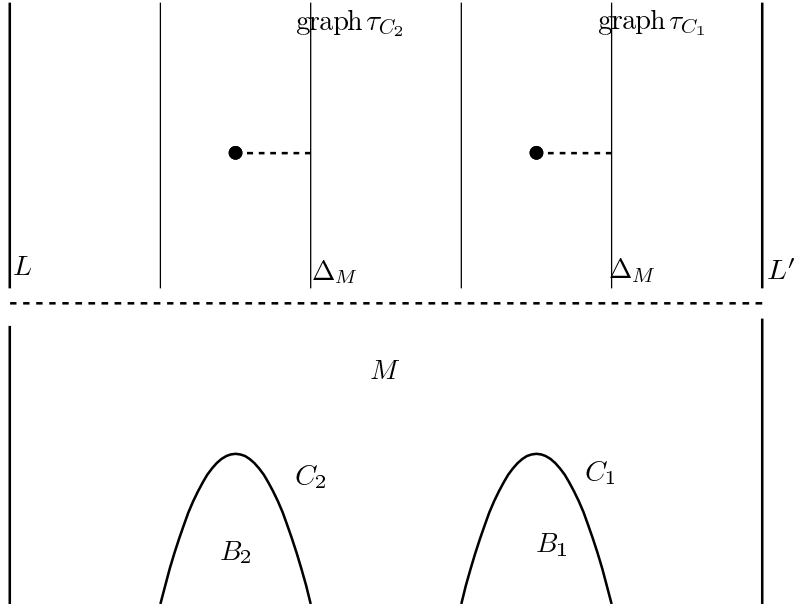}
\caption{The family of quilted fibrations $\Q_{\mathbf{1},\infty}$ (top) which is  used in conjunction with the family $\Q_{I,\mathbf{1}}$ from section \ref{familyofquilts} (bottom) to define the map $f$. The $y$ coordinates of the critical values of the fibration (which are denoted by dotes) can vary but they are bounded above.}\label{qisoquilts}
\end{figure}

%For $\Delta$ its simple to describe. For each $C^{t}C$ it is given by a one parameter family. 
Now we describe the map $f$. The map $g$ has a dual description (as in Remark \ref{negtri}).
Let  $\mathbf{1}=(1,1,\ldots,1)\in \{0,1\}^{n_{+}}$. 
%be such that its $i$'th entry is equal to $1$ if $\epsilon_{i}=+1$ and is zero otherwise. 
Consider the family of quilted fibrations as in Figure \ref{qisoquilts} in which the $y$ coordinates of the critical values are arbitrary but  bounded above.
Denote this family by $\Q_{\mathbf{1},\mathbf{\infty}}$.
The quilt family that gives the map $f_{I}: CF^{+}_{I}\to CF$ is obtained by putting the quilts in the family 
$\Q_{\mathbf{1},\mathbf{\infty}}$ on top of the quilts in the family $\Q_{I,\mathbf{1}}$ from section \ref{familyofquilts} (but now with $n_{+}$ in place of $N$). 
To show that  $f=\sum_{I\in \{0,1\}^{n_{+}}} f_{I}$ is a chain map, since bubbling is ruled out by the Admissibility Condition,
one inspects the limits of the family of quilts used. As in the case $n_{+}=1$ (section \ref{wwtrisection}), when the  $y$ coordinate of a critical value goes to $-\infty$, the resulting map is zero. The remaining boundary components give the equation 
$\sum_{ J\leq I} f_{J}\circ \mu_{I,J}=D \circ f_{I}$ where $D$ is the differential on $CF$. Note that for $n_{+}=1$, $f$ is the map $(h,k)$ from section \ref{wwtrisection}.

%We can describe $f$ as a composition of a bunch of maps.

\bp{(of Theorem \ref{thequasiiso})}
We only prove that $f$ is a quasi-isomorphism by induction on $n_{+}$. The proof for $g$ is similar. %  $g$ is identity on the positive twists and completes the triangle for negative ones. So, it gives a quasi-isomorphism. 
By Theorem \ref{WWtrithm} there is a chain map $(h,k)$ from the cone of
$\mu_{(0),(1)}: CF_{0}\to CF_{1}$ to $CF$  where 
$$CF_{0}=CF(L,C_{N}^{t}C_{N},\, \tau_{C_{N-1}}\circ\cdots\circ \tau_{C_{1}} (L) )$$
and $CF_{1}=CF(L,\Delta_{M},\, \tau_{C_{N-1}}\circ\cdots\circ \tau_{C_{1}} (L) )$ and $(h,k)$ induces an isomorphism on cohomology. 
We use the shorthand notation $iI=(i,I_{1},\ldots, I_{N-1})$ for $i=0,1$ and set $CF_{i\oplus}=\oplus_{I\in \{0,1\}^{N-1}} CF_{iI}$.
By  the induction hypothesis we have  quasi-isomorphisms
$f_{0}: CF_{0\oplus}\to CF_{0}$ and $f_{1}:CF_{1\oplus} \to CF_{1}$. 
%This together with $\mu_{(0),(1)}$ gives us a morphism in the derived category from $CF_{0\oplus} \to CF_{1\oplus}$. This morphism can  be realized by a chain map if
 We construct a  chain map $\nu: CF_{0\oplus}\to CF_{1\oplus}$ which makes the following diagram commutative up to homotopy. 
\begin{equation}\label{commdiag}
\xymatrix{
CF_{0} \ar[r]^{\mu_{(0),(1)}} & CF_{1} \\
CF_{0\oplus}\ar[r]^{\nu} \ar[u]^{f_{0}} & CF_{1\oplus}\ar[u]^{f_{0}}
}
\end{equation}

Set $\nu_{I,J}=\mu_{0I,1J}$.
%Define $\nu_{II}: CF_{0I}\to CF_{1I}$ to be $\mu_{0I,1I}$.
%So we now need to define the rest of the maps $\nu_{IJ}$ which make $\nu$ into a chain map. 
%
%Considering $Cone(\mu_{iI}: CF_{iI}\to CF_{iJ})$ for $i=0,1$, the 
% Lemma \ref{homollem} implies that $\nu_{IJ}$ must be a homotopy between $\mu_{1I,1J}\circ\nu_{II}-\nu_{JJ}\circ \mu_{0I,0J}$ and zero. We note that $\mu_{0I,1J}$ is indeed such a homotopy.
%(For when the $C_{i}$ are Lagrangian spheres this is the \aifty associativity relation.)
%
%Moreover we are free to choose such a homotopy because two such homotopies $\nu_{IJ}$ and $\nu_{IJ}'$ result in chain homotopic maps $\nu$, $\nu'$???
%
That $\nu$ is a chain map is a consequence of Lemma \ref{d2=0}. 
More precisely if  $D_{i\oplus}$ denotes the differential on $CF_{i\oplus}$ then
 $\nu D_{0\oplus}- D_{1\oplus}\nu= D^{2}- D_{0\oplus}^{2}+ D_{1\oplus}^{2}$.
Because the maps $f_{0}, f_{1}$ are given by counting quilts, an argument similar to that of Lemma \ref{d2=0} shows that each $\nu$ makes the  diagram \eqref{commdiag} commutative up to homotopy. This is again a special case of the master equation for family quilt invariants. 
 Therefore $(f_{0},f_{1})$ gives a quasi-isomorphism from $Cone(\nu)= CF_{\oplus}$ to $Cone(\mu_{(0),(1)})$.
Now it is easy to see that the composition $(h,k)\circ (f_{0},f_{1})$ is given by counting the sections of the a family of  quilted fibrations isotopic to the one we used to define $f$. Therefore they induce the same map on homology.
\ep

\begin{remark}
The question of whether the spectral sequence of Theorem \ref{mainthm} collapses at the $E_{2}$ level is related to the formality of the (generalized) Fukaya category of $M$.
Since the $E_{1}$ page is always doubly graded, formality would follow from the existence of a second grading on the Fukaya category which is preserved by the higher composition maps. One can expect that such extra grading should come from extra geometric structures on the manifold e.g. a hyperk\"ahler structure or a circle action.

Note however that the spectral sequence will always collapse at the $E_{N+1}$ term because of the finiteness of the filtration.
\end{remark}

\begin{figure}
\includegraphics[width=0.8\textwidth]{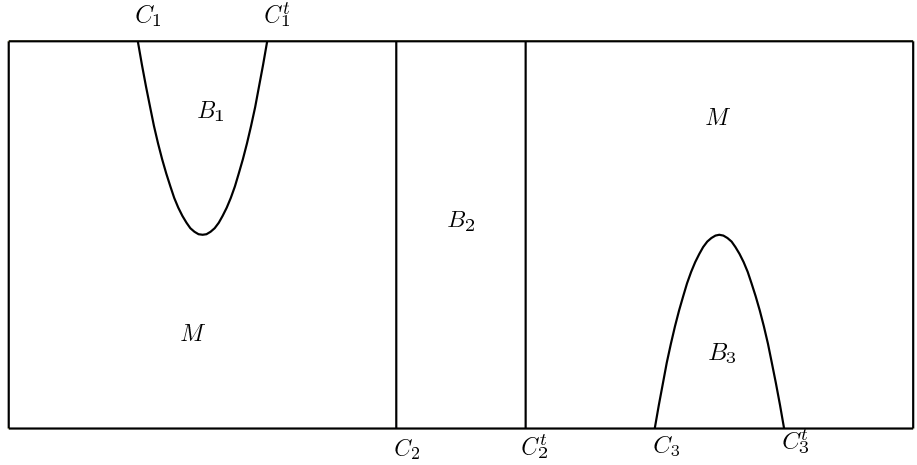}
\caption{An example of the quilts used in constructing the hypercube in $\F(M,M)$. The boundary rectangle represents a cylindrical end.}\label{openquilts}
\end{figure}

For $I\in \{0,1\}^{N}$ denote $\mathbf{C}^{I}= (C_{N}^{I_{N}},\cdots, C_{1}^{I_{1}})$ which is a Lagrangian correspondence from $M$ to itself. Assuming the Admissibility Condition is satisfied by the $C_{i}$, let 
$\tilde\mu_{I,J}\in CF(\mathbf{C}^{I},\mathbf{C}^{J})$  be given by
the count of quilted disks with cylindrical ends which are obtained from quilt $Q_{I,J}$  as in Figure \ref{openquilts}.

\bpr\label{opencube}
Under the same assumptions as in Theorem \ref{thequasiiso}, the Lagrangian correspondence
\begin{equation}
\on{graph} (\tau^{\eps_{N}}_{C_{N}}\circ\cdots \circ \tau^{\eps_{1}}_{C_{1}})
\end{equation}
 is isomorphic, in $D\F^{\#}(M,M)$, to  $(\sum_{I} \mathbf{C}^{I}, \sum_{I\leq J} \tilde{\mu}_{I,J})$.

\epr
Proof is the same as that of Theorem \ref{thequasiiso} but with quilts of the form in Figure \ref{openquilts}.

%\section{Properties of the spectral sequence}

\subsection{Naturality of the spectral sequence}
Since the quasi-isomorphism $\phi$ of Theorem \ref{thequasiiso} is given by quilt maps, it satisfies a form of naturality under equivalence of Lagrangian correspondences as follows.  First recall that if $L_{0}\subset M_{0}\times M_{1}$ and $L_{1}\subset M_{1}\times M_{2}$ are two Lagrangians correspondences then their \emph{composition} $L_{1}\circ L_{0}\subset M_{0}\times M_{2}$ is 
$$\{(m_{0},m_{2}) |\, \exists m_{1}\in M_{1}\, s.t.\, (m_{0},m_{1})\in L_{0}\,\&\, (m_{1},m_{2})\in L_{2}\}.$$ 
It can be equally described as the intersection of $L_{0}\times L_{1}$ with $M_{0}\times\Delta_{M_{1}}\times M_{2}$ in $M_{0}\times M_{1}\times M_{1}\times M_{2}$. We say that the composition of $L_{0}$ and $L_{1}$ is \emph{embedded} if this intersection is transversal and its projection into $M_{0}\times M_{1}$ is an embedding.

Let 
\begin{equation}\label{L}
\mathbf{L}=(L_{1},L_{2},\ldots, L_{n})
\end{equation}
 be a generalized Lagrangian correspondence. Assume further that for some $1\leq k\leq n-1$, the composition $L_{k+1}\circ L_{k}$ is  embedded.
Then $\mathbf{L}$ and
\begin{equation}\label{L'}
 \mathbf{L'}=(L_{1} ,\ldots, L_{k-1} ,L_{k+1}\circ L_{k}, L_{k+2}, \ldots, L_{n})
\end{equation}
 are said to be equivalent.
More generally  $\mathbf{L}$ and $\mathbf{L}'$ are said to be \emph{equivalent in the symplectic category} if they are related by a sequence of such moves.
It is not difficult to see that the generators of $CF(\mathbf{L})$ and $CF(\mathbf{L}')$ are in one to one correspondence.
The functoriality theorem of Wehrheim and Woodward \cite[Thm. 1.0.1]{Functoriality} (together with the discussions of monotonicity in \cite{QuiltedFloer}) %(and generalization in \cite{})
 implies that if $\mathbf{L}$ and $\mathbf{L'}$ satisfy the Admissibility Condition then their Floer homologies are canonically isomorphic. Moreover such isomorphisms are compatible with maps induces by quilts. More precisely we have the following result (which is stated in greater generality in  \cite{WWquilts}).

\bthm[Wehrheim, Woodward \cite{WWquilts}, Thm. 5.1]\label{WWnaturality}
Let $Q$ be a quilt with one incoming and  one outgoing end. Assume that the incoming (resp. outgoing) end is labeled by $\mathbf{L}_{0}$ (resp. $\mathbf{L}_{1}$). Assume that $Q$ has a strip in it whose seams are labeled by correspondences $L_{k}$ and $L_{k+1}$ whose composition is embedded. Let $Q'$ be the quilt obtained from $Q$ by removing the this strip and replacing it with a seam labeled by $L_{k+1}\circ L_{k}$. Let $\mathbf{L}'_{0}, \mathbf L_{1}'$ be the labellings of the incoming and outgoing ends of $Q'$ respectively.  Then one has a commutative diagram
\begin{equation}
\xymatrix{
CF(\mathbf{L}_{1}) \ar[r]  & CF(\mathbf{L}'_{1})\\
CF(\mathbf{L}_{0}) \ar[r] \ar[u]^{\Phi(Q)}  & CF(\mathbf{L}'_{0}) \ar[u]_{\Phi(Q')}
}
\end{equation}
where  the horizontal maps are isomorphisms and vertical ones are quilt maps.
 Furthermore if $\mathbf{L}_{0},\mathbf{L}_{1}$ satisfy the Admissibility Condition \eqref{wwfuk} the horizontal maps  induce isomorphisms on homology.
\ethm

Note that Lemma \ref{nobubble} together with the discussions of monotonicity in \cite{WWquilts} imply that the above theorem holds if $\mathbf{L}_{i}$ satisfy the Admissibility Condition.

\bpr\label{naturalfunct}
Let $\mathbf{L}_{i},\mathbf{L}'_{i}$ be two generalized submanifolds of $M$ which are equivalent in the symplectic category for $i=0,1$ , and $\mathbf{D},\mathbf{D}'$ two equivalent generalized correspondence from $M$ to itself. Finally let $\psi_{0},\psi_{1}$ be two symplectomorphisms of $M$ given by compositions of fibered Dehn twists along spheric coisotropic submanifolds of $M$. 
If all these satisfy the Admissibility Condition then we have a commutative diagram %of chain maps
\begin{equation}
\xymatrix{
CF(\mathbf{L}_{0},\psi_{0},\mathbf{D},\psi_{1},\mathbf{L}_{1})\ar[r] & CF(\mathbf{L}'_{0},\psi_{0},\mathbf{D}',\psi_{1},\mathbf{L}'_{1})\\
CF_{\oplus}(\mathbf{L}_{0},\psi_{0},\mathbf{D},\psi_{1},\mathbf{L}_{1})\ar[r] \ar[u] & CF_{\oplus}(\mathbf{L}'_{0},\psi_{0},\mathbf{D}',\psi_{1},\mathbf{L}'_{1})\ar[u]
}
\end{equation}
where the horizontal maps are natural isomorphisms  which induce isomorphisms on homology.
 The vertical ones are isomorphisms (in the derived category) given  by Theorem \ref{thequasiiso}. %$HF_{\oplus}$ represents.
\epr
The proof is an application of Theorem \ref{WWnaturality}.
 %\red{Monotonicity, Yanki's proof}

%\grn{The effect of canceling a correspondence on the hypercube}

%\begin{equation}
%CF(L, \tau_{\phi}, D^{t},D,\tau_{\psi}, L') \to %CF(L,\tau_{\phi},\tau_{\psi},L')
%\end{equation}

\section{Applications}

\subsection{ Symplectic Khovanov homology}\label{khsexample}
\renewcommand{\L}{\mathcal{L}}
Symplectic Khovanov homology is an invariant of links introduced by Seidel and Smith \cite{SS}. It is expected to be equivalent to Khovanov homology. We use the same notation as in section \ref{ssex}. Let $z_{j}=(j,0)$ for $j=0,\ldots, 2m$ be points in the plane. They give rise to a point $\nu\in Conf_{2m}$. Let $\delta_{i}:[0,1)\to Conf_{2m}$ be a curve such that $\delta(0)=\nu$ and as $t\to 1$, $z_{i}$ and $z_{i+1}$ merge  (but the other points remain fixed). As mentioned in section \ref{ssex}, this gives us locally defined spheric coisotropic submanifolds $L_{i}=L_{\delta_{i}}\subset \Y_{m,\nu}$ which fiber over a compact subset $U\subset \Y_{m-1,\mu}$ where $\mu=\nu\bs \{z_{i},z_{i+1}\}$.  By composing the correspondences $L_{2i-1}$ for $i=1,\ldots, m$ we get a Lagrangian submanifold $\L\subset \Y_{m,\nu}$ which does not depend on the choice of the open sets $U$ (because it is compact). We set the compact subset $U\subset \Y_{m-1,\mu}$ %(used in the definition of the $L_{\delta}$) 
to contain the projection  of a \nbhd of $\L$. Even though fibered Dehn twists along the $L_{j}$ are defined locally, since $\L$ is compact and the Dehn twists are compatible (up to Hamiltonian isotopy) under enlargement of $U$, the image of $\L$ under a composition of such Dehn twists is a well-defined Lagrangian submanifold up to Hamiltonian isotopy.

 %Also because the first Chern class of $\Y_{m,\nu}$ vanishes
%
Following the (slight) reformulation in \cite{RR1} let a link $K$ be given as the plat closure of a braid $\beta\in Br_{2m}$ which in turn is given by a braid word $\sig^{\eps_{1}}_{k_{1}}\sig^{\epsilon_{2}}_{k_{2}}\cdots \sig^{\epsilon_{N}}_{k_{N}}$. Let $w$ denote the writhe of this braid. One has
\begin{equation}
Kh_{s}(K)=HF(\L,\tau^{\epsilon_{1}}_{L_{k_{1}}}\circ\cdots\circ \tau^{\epsilon_{N}}_{L_{k_{N}}} (\L))[-m-w].
\end{equation}

\begin{remark}
Seidel and Smith define their invariant using ``rescaled''  monodromy maps of the fibration $\chi$. The equivalence with the above formulation follows from a theorem of Perutz \cite[Theorem 2.19]{PerutzI}, which states that the monodromy of a normally K\"ahler Lefschetz-Bott fibration around a critical value is Hamiltonian isotopic to fibered Dehn twist along the vanishing cycle of the fibration, and the fact that the Lagrangian $\L$ is compact.  
\end{remark}

% The spectral sequence in Section \ref{spectralseq} applied to ${Kh}_s$ gives a spectral sequence from the hyper-cube of resolution for ${Kh}$ to ${Kh}_s$
The Lagrangian correspondences $L_{i}$ and the Lagrangian submanifold $\L$ are exact because the symplectic form on $\Y_{m,\nu}$ is exact and the fibers of the fibration $L_{i}\to \Y_{m-1,\bar{\nu}}$ are simply connected. 
It was shown in \cite[Section 3]{RR1} that one can make the invariant under the Liouville flow outside a compact subset (which we take to be $U$) without affecting Floer homology.
Therefore  they are admissible of type \eqref{myfuk} and (with Lemma \ref{mycompact} in mind) one can apply Theorem \ref{mainthm}. 
The maps $\mu_{I,J}$ between the adjacent vertices of the hypercube, which are given by counting quilted triangles, where used in \cite{RR2} to define homomorphism $Kh_{s}(S)$ corresponding to elementary cobordisms $S$ between links (and more generally tangles).
 To identify the the $E_{1}$ page we use the following result which is proved in Theorem 5.6 in \cite{RR2}. (See also Prop. 5.12 therein.)

\bpr%[\cite{RR2}, Thm. 5.5]
Let $K,K'$ be two flat unlinks in the plane, related by an elementary cobordism $S$. Then, with coefficients in $\Z/2$,   one has a commutative diagram
\begin{equation}
\xymatrix{
Kh_{s}(K)\ar[r]^{Kh_{s}(S)} \ar[d] & Kh_{s}(K') \ar[d]\\
Kh(K) \ar[r]^{Kh(S)} & Kh(K')
}
\end{equation}
where the vertical arrows are isomorphisms.
\epr

Therefore we obtain the following.
\bpr
For each link $K$ there is a spectral sequence whose $E_{2}$ term is the Khovanov homology with coefficients in $\Z/2$ of $K$ and converges to $Kh_s(K)$. 
\epr

%$\mathcal{H}_m : = \bigoplus \mathcal{H}_{\mathcal{Y}_m} = %\underset{a,b}{\bigoplus} CF(L_a, L_b)$

%Taking cohomology $\mathcal{H}_m \longrightarrow$Khovanov's arc rings $H_m$

%For ${Kh}_s$ and for a link $K = T_k T_{k - 1} \cdots T_0$ we get
%{ \begin{equation}
%  {Kh}_s (K) \cong {Kh}_s (T_k)
%  \underset{\mathcal{H}_{m_k}}{\otimes_{}} \cdots
%  \underset{\mathcal{H}_{m_1}}{\otimes} {Kh}_s (T_0)
%\end{equation}
%}

%\grn{Use abouzaid dg to show that spectral sequence vanishes for...}
%At last in this subsection, we use a method due to  Abouzaid \cite{} to reobtain  a result of Waldron \cite{} that the \aifty ring $H_{2}$ is formal. This is motivated by a result of Lekili and Perutz \cite{} on the \aifty structure of the Fukaya category of the torus which uses the same method. However the computation in our case is simpler than that of \cite{} (because of formality).  Waldron uses Kuranishi homology to obtain this result but we do not. 

\subsection{\HG homology} \label{hgproof}
In this section we prove Theorem \ref{HGss}.
Let $\Sigma$ be a surface of genus $g$ and equip $Sym^{g} \Sigma$ with the K\"ahler form in the cohomology class $P_{\Sigma}$ from Example \ref{perutzex}.
For the coisotropic submanifolds $V_{\gamma}\subset Sym^{g} \Sigma$ the inclusion map is injective on the fundamental groups and therefore $\pi_{2}(Sym^{g} \Sigma, V_{\gamma})=0$.
Since (with notation from Example \ref{perutzex}) $c_{1}(Sym^{k}\Sigma)=(k+1-g)\eta-\theta$ and the integral of  $\theta$ over spheres is zero, the above K\"ahler form makes $Sym^{g} \Sigma$ into a spherically monotone symplectic manifold.
 The submanifolds $V_{\gamma}$ are also balanced; the needed line bundle is given by the anticanonical bundle of $Sym^{g}\Sigma$. See for example Section 6.1 in \cite{DenisFuk-HF}. Therefore they are admissible of type \eqref{balancedfuk}.

%The exact triangle of Theorem \ref{WWtrithm} is stated for coisotropic submanifolds of codimension at least two however this condition is to   
%{WW need codim>1}

Let $M$ be a 3-manifold given by gluing a genus $g$ handlebody $H$ to another such handlebody $H'$ by $ \phi$ where $\phi$ is an element of the mapping class group of $\Sigma=\Sigma_{g}=\partial H$.
 %and $\iota$ is the orientation reversing homeomorphism of $\Sigma$ that sends the meridians 
Let $H$ and $H'$ be given respectively  by attaching disks to circles $\alpha_{1},\ldots, \alpha_{g}$ %in Figure \ref{abcurves}  and $H'$ by arbitrary curves  
  and $\beta_{1},\ldots,\beta_{g}$  in $\Sigma$. It is easy to see that $(\Sigma,{\alpha'},{\beta})$, where 
${\alpha}'=(\phi(\alpha_{1}),\ldots,\phi(\alpha_{g}))$ and 
${\beta}=(\beta_{1},\ldots,\beta_{g})$, form a Heegaard diagram for $M$.  Let $T_{\alpha}=\alpha_{1}\times \alpha_{2}\times \cdots\times \alpha_{g}$ and
 $T_{\beta}=\beta_{1}\times \beta_{2}\times \cdots \times \beta_{g}$.%  and let $\phi_{*}$ denote the map induced on $Sym^{g}\Sigma$ by $\phi$.  

\bl\label{callahantype}
If $\gamma$ is an embedded circle in $\Sigma$ and the coefficient of $\theta$ in $P_{\Sigma}$ is positive then
$\tau_{V_{\gamma}} (T_{\alpha})$ is Hamiltonian isotopic to 
$T'_{\alpha}=\tau_{\gamma} \alpha_{1}\times \tau_{\gamma} \alpha_{2}\times \cdots \times \tau_{\gamma}\alpha_{g}$.

\el
\bp
By a result of Lekili  \cite[3.4.1]{lekili:thesis}, $T_{\alpha}$ is Hamiltonian isotopic to $V_{\alpha_{1}}\circ\cdots \circ V_{\alpha_{g}}$. One way to prove the lemma is to note that by the naturality of vanishing cycle construction (and the isotopy of monodromy with fibered Dehn twist along the vanishing cycle) we have
\bq
\tau_{V_{\gamma}}( V_{\alpha})\cong V_{\tau_{\gamma} \alpha}. %\circ \on{grapg} \tau_{V_{\gamma}}.
\eq

Alternatively we can consider a Lefschetz fibration $p$ over the unit disk in which the curve $\gamma$ gets pinched to a point. Let $\pi$ denote  the relative Hilbert scheme of the fibration $p$. 
Let $\alpha^{t}_{i}\subset p^{-1}(e^{ it})$ for $1\leq i\leq g$ be a smooth family so that $\alpha^{0}_{i}=\alpha_{i}$ and $\alpha^{2\pi}_{i}= \tau_{\gamma}\alpha_{i}$.  
We can construct a Heegaard torus 
$T_{\alpha,t}:=V_{ \alpha^{t}_{1}}\circ \cdots\circ V_{\alpha^{t}_{g}}\subset \pi^{-1}(e^{it})$.
By pulling back these tori into $\pi^{-1}(1)$ using parallel transport maps, we get a Lagrangian isotopy between $T_{\alpha}$ and $\kappa^{-1}(T'_{\alpha})$
where $\kappa$ is the monodromy of $\pi$ which by \cite[Thm 2.19]{PerutzI} is Hamiltonian isotopic to fibered Dehn twist along $V_{\gamma}$.
Since an isotopy through balanced Lagrangians is always given by a Hamiltonian isotopy, the result follows. 

% To see this note that if $\psi$ and $W$ are respectively a symplectomorphism and a spheric coisotropic submanifold of a symplectic manifold $M$ then $\tau_{\phi(W)}\cong\psi\circ \tau_{W} \circ \psi^{-1}$. Letting $\psi=\tau_{V_{\gamma}}$ and $W=V_{\alpha}$ the (unique up to isomorphism) Lefschetz-Bott fibration over the disk whose monodromy is given by $\tau_{\psi(W)}$ has $\psi(W)$ as vanishing cycle. Also if $\kappa$ is the homomorphism introduces  in section 1.1 of \cite{PerutzII} then by \cite[Thm. 1.1]{PerutzII}, $\psi^{-1}\circ \tau_{W} \circ \psi\cong \kappa(\tau_{\gamma} \tau_{\alpha}\tau_{\gamma}^{-1})$. The latter is by definition the symplectic monodromy of the relative Hilbert scheme of a Lefschetz fibration $P$ whose symplectic monodromy is $\tau_{\gamma}\tau_{\alpha}\tau_{\gamma}^{-1}=\tau_{\tau_{\gamma} \alpha}$. The vanishing cycle of this relative Hilbert scheme is, again by definition, the same as $V_{\tau_{\gamma}\alpha }$.  
\ep

Now let $\phi=\tau_{\gamma_{k}}\circ \cdots \circ\tau_{\gamma_1}$ be an expression of $\phi$ as a composition of (classical) Dehn twists along a number of curves $\gamma_{1},\ldots,\gamma_{k}$ in $\Sigma$. Set $\phi_{*}=\tau_{V_{\gamma_{k}}}\circ \cdots \circ \tau_{V_{\gamma_{1}}}$. It follows from Lemma \ref{callahantype} that $T_{\alpha'}$ is Hamiltonian isotopic to %$\tau_{V_{\gamma_{k}}}\circ \cdots \circ \tau_{V_{\gamma_{1}}} (T_{\alpha})$.
$\phi_{*}(T_{\alpha})$. 
As in \cite{perutz:hamilslides} the Heegaard tori are Lagrangian for our choice of the symplectic form.
Let $h$ be a Hamiltonian diffeomorphism which makes the Heegaard tori admissible. Therefore  
the hat version of the \HG homology of $M$ is given by
\begin{equation}
\widehat{HF}(M)= HF(T_{\beta},h\circ\phi_{*}(T_{\alpha}))
\end{equation}
where the Floer homology is taken in $Sym^{g}(\Sigma)\bs (z\times Sym^{g-1}\Sigma)$.
%\red{Welldefined? Must make sure the diagram is admissibile; for example if $\phi$ sends $\alpha_{i}$ to itself.}
%
Let $n_{z}$ denote the intersection number with the hypersurface $z\times Sym^{g-1}\Sigma$.
Note that a sequence of holomorphic curves with $n_{z}=0$ will not converge to a curve with $n_{z}\neq 0$ because of the additivity of $n_{z}$ and the fact that $n_{z}\geq 0$ for holomorphic curves.
 Theorem \ref{HGss} now follows from Theorem \ref{thequasiiso}.

 \subsection{ The spectral sequence of a branched double cover}\label{osproof}
\begin{figure}[ht]
\includegraphics{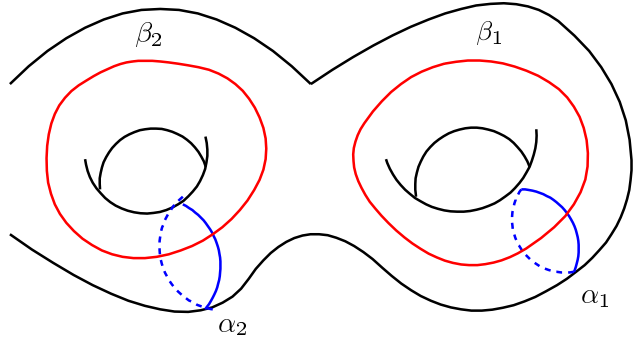}
\caption{Special $\alpha$ and $\beta$ curves}\label{abcurves}
\end{figure}

In this section we show that %Ozsvath and Szabo's 
the spectral sequence for branched double covers of links is a special case of the spectral sequence of Theorem \ref{HGss}. Let $K$ be a link in $S^{3}$ given as the plat closure of a braid $b$ on ${2m}$ strands. 
We add two auxiliary strands to $b$ which are not linked with each other or with other strands. 
% Without thisthings won't work. 
 We  denote the resulting $2m+2$ braid by $b'$.
Let $B\subset S^{3}$ be a ball \st $B\cap K=b'$ and $S=\partial B$ intersects $K$ transversely in $4m+4$ points.
The branched double cover of $S^{3}\bs B$ consists of two handlebodies whose boundary is a surface $\Sigma$ of genus $m$. It is given  by attaching disks to the $\beta$ curves in Figure \ref{abcurves}.

The double cover of $B$ branched over $b'$ is the mapping cylinder of an element $\phi$ of the mapping class group of $\Sigma$ which is given as follows. Let $b=\sig^{\epsilon_{N}}_{k_{N}}\cdots \sig^{\epsilon_{1}}_{k_{1}}$ be an expression of $b$ in terms of  braid generators.
Let $\beta_{1},\ldots, \beta_{m}$ be as in Figure \ref{abcurves} and $\gamma_{1},\ldots, \gamma_{m-1}$ be curves in $\Sigma$ \st $\gamma_{i}$ meets each one of $\beta_{i}$ and $\beta_{i+1}$ in exactly one point and does not intersect other $\beta$ curves. 
Let  $\delta_{2i-1}$ be equal to $\beta_{i}$ for $i=1,\ldots, m$ and  $\delta_{2i}$ equal to ${\gamma_{i}}$ for $i=1,\ldots, m-1$.
 Then $$\phi=\tau^{-\epsilon_{N}}_{\delta_{k_{N}}}\circ\cdots \circ\tau^{-\epsilon_{1}}_{\delta_{k_{1}}}.$$
 This is because the conventions for positive braids and positive Dehn twists are opposite of each other.

%$\tau^{\epsilon_{k}}_{\beta_{n_{k}}}\cdots \tau^{\epsilon_{1}}_{\beta_{n_{1}}}$ where the $\beta_{i}$ are as in Figure \ref{abcurves}.

%Now it follows from Theorem \ref{HGss} that there is a spectral sequence starting from \red{$\oplus_{I} HF(T_{\alpha,} C^{I},T_{\alpha})$} and converging to $\HFh(\Sigma(K))$.

\bl[Lekili, Perutz]\label{commrel}
%Let $\alpha_{1},\ldots,\alpha_{m}$ be the collection of $\alpha$ circles on $\Sigma_{m}$.
 The Perutz correspondences  associated to the curves $\beta_{i}, \gamma_{i}$ satisfy the following relations where $\cong$ denotes Hamiltonian isotopy.
\begin{enumerate}
\item\label{tripint} $V_{\beta_{i}}\circ V_{\gamma_{i+1}}^{t}\cong \Delta_{Sym^{m-1}\Sigma_{m-1}}$
\item\label{compprod} $V_{\beta_{1}}\circ\cdots\circ V_{\beta_{m}}\cong \beta_{1}\times \cdots\times \beta_{m}$ 
\end{enumerate}
\el
\
Both statements follow from work in progress of Lekili and Perutz
\cite{lekili-p} and \eqref{compprod} is a special case of a theorem
proved by Lekili  in his thesis \cite[3.4.1]{lekili:thesis} for a general set of nonintersecting curves. 
Note that (i) is obvious in case $m=1$. The proof for the general case uses a delicate degeneration argument.
%However note that because of the special form of the curves $\beta_{i}$, \eqref{compprod} follows, using induction, from the fact that the bundle $V_{\gamma}\to Sym^{m-1} \Sigma_{m-1}$ 
%(as in the proof of Lemma \ref{callahantype}) is trivial as a holomorphic bundle \cite[Prop. 3.13 (5)]{PerutzI}. It is also enough to prove \eqref{tripint} for the case $m=1$ which becomes trivial.

For $I\in\{0,1\}^{I}$ let $\mathbf{V}^{I}$ denote the generalized Lagrangian correspondence $(T_{\alpha}, V_{\delta_{N}}^{I_{N}}, \ldots V_{\delta_{1}}^{I_{1}} ,T_{\alpha})$.
The above lemma together with the obvious relation $\tau_{V_{\delta}} \circ \tau_{V_{\delta}}^{-1}=id$
are enough to  conclude that each $\mathbf{V}^{I}$ is equivalent in the symplectic category to a generalized correspondence of the form 
\begin{equation}\label{endcorr}
(\beta_{1}\times \cdots \times \beta_{k}, \beta_{1}\times \cdots \times \beta_{k})
\end{equation} 
in $Sym^{k} \Sigma_{k}$ with $k\leq m$. Basically the relations of Lemma \ref{commrel} imply that the correspondences $V_{\delta_{i}}$ satisfy the same commutation relations as the flat tangles they are assigned to.
Let $A$ be the ungraded Khovanov's algebra $H^{1}$ over $\Z/2$ i.e.~it is generated over $\Z/2$ by $1,X$ with relations $X^{2}=0,\,
11=1,\, 1X=X$ and with no grading.
After Hamiltonian isotoping the above correspondence to make it balanced (which is the same as admissibility of the corresponding Heegaard diagram in this case), the Lagrangian Floer homology of \eqref{endcorr} is, as a vector space, isomorphic to  $A^{\tens k}$. This is because $\beta_{i}$ and a Hamiltonian isotoped copy of it intersect at two points which are both cocycles.

Let $Kh_{k,k+1}: A^{\tens k} \to A^{\tens k+1}$  and $Kh_{k,k-1}: A^{\tens k} \to A^{\tens k-1}$ be given by Khovanov's TQFT. 
More specifically $Kh_{2,1}$ is the multiplication map on $A$ and $Kh_{1,2}$ sends $1$ to $1\tens X+X\tens 1$ and $X$ to $X\tens X$.
 Our notation is misleading because it does not specify on which factor of $A^{\tens k}$ the maps are acting but this will be clear from the context.

\begin{figure}[ht]
\includegraphics{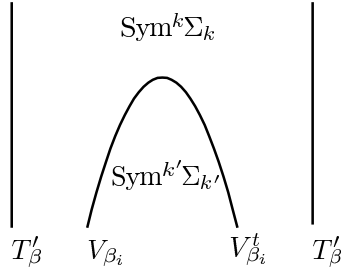}
\caption{The quilted triangle $Q'_{I,J}$ used in the proof of Prop. \ref{hgcomdiag}. Here $k'=k-1$. }\label{quilttrisym}
\end{figure}

\newcommand{\bV}{\mathbf{V}}
\bpr\label{hgcomdiag}
If $J$ is an immediate successor for $I$ then we have a commutative diagram where the vertical arrows are isomorphisms and $k'=k\pm 1$.
\begin{equation}
\xymatrix{
HF(\mathbf{V}^{I}) \ar[r]^{\mu_{I,J}} \ar[d] & HF(\bV^{J})\ar[d]\\
A^{\otimes k} \ar[r]^{Kh_{k,k'}} & A^{\otimes k'}
}
\end{equation}
\epr
\bp
Using the relations of Lemma \ref{commrel} (together with $\tau_{V}\tau_{V}^{-1}=id$) 
the Lagrangian correspondences $\bV_{I}$ and $\bV_{J}$ become equivalent to correspondences either of the form $\mathbf{W}=(T'_{\beta},V_{\beta_{j}},V_{\beta_{j}}^{t}, T'_{\beta})$ for some $j$ or $\mathbf{U}=(T'_{\beta},T'_{\beta})$ where $T'_{\beta}={\beta_{1}}\times \cdots\times \beta_{k}$ for some $k\leq m$.
The set of equivalences that give $\mathbf{U}$ and $\mathbf{W}$ from $\bV^{I}$ and $\bV^{J}$   correspond to a sequence of strip shrinking (and/or unshrinking) in the quilt $Q_{I,J}$ which results in  a quilted pair of pants $Q'_{I,J}$ which is labelled by $\mathbf{U}$ and $\mathbf{W}$ as in Figure \ref{quilttrisym}. We apply Theorem \ref{WWnaturality} to $Q_{I,J}$ and $Q'_{I,J}$ which 
reduces the problem to computing the maps induced by $Q'_{I,J}$.
%
%Now we use Theorem \ref{myversiontri} which tells us that the moduli space of holomorphic maps of $Q'_{I,J}$ is cobordant to  the moduli space of \pse triangles $u$ in $Sym^{k}\Sigma_{k}$  with boundary conditions given by Lagrangians $T'_{\beta}, h(T'_{\beta})$ and the coisotropic $V_{\beta_{j}}$ where $h$ is a Hamiltonian isotopy.

%One can use almost complex structures on $Sym^{k}\Sigma_{k}$ induced by such structures on $\Sigma_{k}$ to achieve transversality for maps $u$. (An argument for a similar statement is given in \cite{lipshitz:cylindrical}.)

A  complex structure (instead of a family) is enough to achieve transversality for transversely intersecting curves on surfaces as in \cite[13b]{Seidelbook}.  
We show that the moduli space of \pse maps of $Q'_{I,J}$, given by complex structures on $Sym^{k}\Sigma_{k}$ and $Sym^{k'}\Sigma_{k'}$ induced from those on $\Sigma_{k} $ $\Sigma_{k'}$, can be described in terms of triangles in $\Sigma_{k}$ and therefore one can deduce the transversality of the former moduli space.
% $Sym^{k}\Sigma_{k}$ and $Sym^{k'}\Sigma_{k'}$ induced from regular strcutures on $\Sigma_{k}$ and $\Sigma_{k}'$, the relevant moduli spaces are smooth. 
%
With the same notation as in Prop.~\ref{myversiontri}, the moduli space of holomorphic maps of $Q'_{I,J}$ is the same as the moduli space of \pse triangles $u$ in $Sym^{k}\Sigma_{k}$  with boundary conditions given by Lagrangians $T'_{\beta}, h(T'_{\beta})$ ($h$ is a Hamiltonian isotopy) and the coisotropic $V_{\beta_{j}}$ for which $\pi\circ u|_{\gamma}$ can be filled with a \pse disk $v$ in $Sym^{k'}\Sigma_{k'}$.

With the choice of complex structures as above, $u$ is given by a map $u'$ of a suitable branched cover of the disk into $\Sigma_{k}$ itself. To avoid the basepoint this map $u'$ has to consist of $k$ holomorphic maps from the disk into $\Sigma_{k}$. 
% This is because of the special form of the Heegaard diags.
Therefore $u$ is determined by $k$ holomorphic maps $u_{1},\ldots, u_{k}$ from the disk into $\Sigma_{k}$. To avoid the basepoint, (after possible renumbering)  $u_{i}$ has to send the boundary components of $Q'_{I,J}$ to $\beta_{i}$ and a Hamiltonian isotopic copy of it. For the same reason there is a unique $1\leq j\leq k$ such that  if $i\neq j$ the image of $u_{i}$ does not intersect $\beta_{j}$ and therefore, 
%using a similar argument to that of Lemma \ref{callahantype},
 the (interior) seam condition $V_{\beta_{j}}$ does not impose any further restriction on $u_{i}$
 \footnote{This is because if $z_{1},\ldots, z_{j-1},z_{j+1},\ldots z_{k}$ are points on $\Sigma_{k}$ away from $\beta_{j}$ and $z_{j}\in \beta_{j}$ then $(z_{1},\ldots,z_{k})\in V_{\beta_{j}}$. 
% This is because $V_{\beta}$ is a bundle over $\Sigma_{k-1}$.
}
 but it implies that $u_{j}$ sends the seam to a third isotopic copy of $\beta_{j}$.

\begin{figure}
\centering
\hspace*{\fill}
{\includegraphics{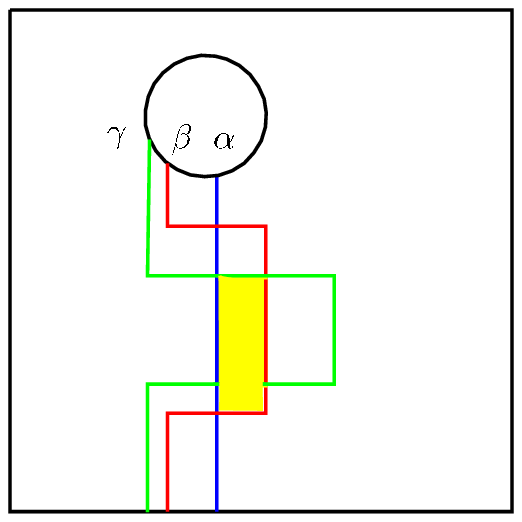}}
\hfill
{\includegraphics[scale=0.87]{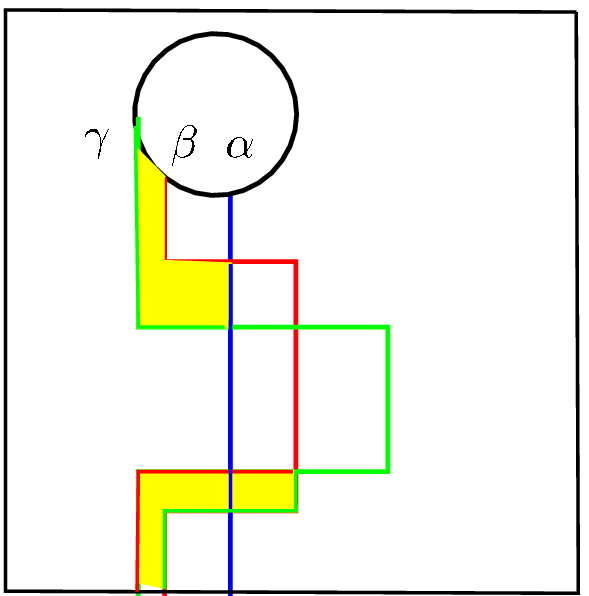}}
\hspace*{\fill}
\caption{The two holomorphic triangles (with one incoming and two outgoing ends) which tell us that $\theta^{+}$ is sent to $\theta^{+}\tens \theta^{-}+\theta^{-}\tens \theta^{+}$ by $\mu_{(0),(1)}$ (for a negative twist). }
\label{heegtries}
\end{figure}

Note that for such a $u$ the disk $v$ always exists because $\pi\circ u|_{\gamma}$ consists of $k'$ nullhomotopic circles in $\Sigma_{k'}$. (This in particular implies that each $u_{i}$ for $i\neq j$ is a \pse strip of Maslov index zero and is therefore constant.)
%
% refer to lipshitz
%Therefore \pse representations of $Q'_{I,J}$ are in one to one correspondence with \pse triangles in $\Sigma_{1}$ with boundary on three Hamiltonian isotopic copies of $V_{\beta_{1}}=\beta_{1}$. 
%
%There are two cases depending on whether the mutation from $\bV^{I}$ to $\bV^{J}$ is of the form $V_{\delta_{i}}^{t} V_{\delta_{i}}\to \Delta_{\Sigma_{m}}$ or the other way around.
%In view of the fact that  
Therefore the computation is reduced to the genus one case and there are two cases to consider: $k=2,k'=1$ and $k=1,k'=2$. (Note that the symplectic form induces on a  punctured torus is exact and therefore the balanced assumption on the Lagrangians is reduced to exactness.)
%
%The rest of the proof is similar to that of Thm. 6.3 in \cite{BranchedDouble} but with a minor difference since the \pse triangles we are considering do not exactly give cobordism maps in \HG homology.
For the first case one can see by direct inspection of the genus one Heegaard diagram that the pair of pants map acts on the generators in the same way as the multiplication in $A$.
This can also bee seen from the fact that $\HFh(S^{1}\times S^{2})$ is freely generated,  over $H^{*}(S^{1})$, by one element. The second case is also proved by inspecting the Heegaard diagram. See Figure \ref{heegtries}.

\ep

Now the Theorem \ref{osthm} follows from Theorem \ref{HGss}  using Lemma \ref{hgcomdiag} and Prop. \ref{naturalfunct}. Note that the extra two strands we added to $b$ adds an unlinked unknot to $K$ and we have $\Sigma(K\cup
 \bigcirc )=\Sigma(K)\# S^{1}\times S^{2}$.

\bibliographystyle{plain}
\bibliography{../biblio}
\vspace{1cm}
Laboratoire de Math\'{e}matiques Jean Leray, Universit\'e de Nantes \\  %\\\newline France
\verb"rezazadegan-r@univ-nantes.fr"

\end{document}